\documentclass[12pt]{amsart}

\usepackage{amssymb,amscd,pb-diagram}
\usepackage{fullpage}
\usepackage[all]{xy}

\newcommand{\C}{{\mathbf C}}
\newcommand{\E}{{\mathcal E}}
\newcommand{\F}{{\mathbf F}}

\newcommand{\M}{{\mathfrak M}}

\newcommand{\Q}{{\mathbf Q}}

\newcommand{\Z}{{\mathbf Z}}

\newcommand{\m}{{\mathfrak m}}
\newcommand{\n}{{\mathfrak n}}
\newcommand{\p}{{\mathfrak p}}

\newcommand{\af}{{\alpha\in_\text{af} K}}
\newcommand{\be}{{B_\lambda}}

\newcommand{\ov}[1]{{\overline{{#1}}}}
\newcommand{\nl}{{\mathcal L_n}}

\newcommand{\pf}[1]{{\mathbf P_K^{#1}}}
\newcommand{\oo}{{\mathcal O}}

\newcommand{\po}{{\mathbf P_K^1}}

\newcommand{\pn}{{\mathbf P_K^n}}

\newcommand{\gal}{{ \mbox{\rm Gal}}}

\newcommand{\vsp}{\vspace{8pt}}

\newcommand{\rarr}{{\rightarrow}}

\newtheorem{thm}{Theorem}
\newtheorem{lem}{Lemma}

\newtheorem{prop}{Proposition}

\newtheorem{defn}{Definition}

\theoremstyle{remark}
\newtheorem{remark}{Remark}

\newtheorem*{ack}{Acknowledgments}

\newcount\timehh\newcount\timemm
\timehh=\time 
\divide\timehh by 60 \timemm=\time
\begin{document}

\title[Specializations of one-parameter families of polynomials]
      {Specializations of one-parameter\\families of polynomials}

\author{Farshid Hajir}
\email{hajir@math.umass.edu}
\address{Department of Mathematics \& Statistics, University of Massachusetts.
	Amherst, MA 01003-9318 USA}

\author{Siman Wong}
\email{siman@math.umass.edu}

\address{Department of Mathematics \& Statistics, University of Massachusetts.
	Amherst, MA 01003-9318 USA}





\thanks{Hajir's research is supported in part by NSF Grant No.~0226869}

\begin{abstract}
Let $K$ be a number field, and let
$
\lambda(x,t)\in K[x, t]
$
be irreducible over $K(t)$.
Using algebraic geometry and group theory, we study the set of
$\alpha\in K$ for which the specialized polynomial $\lambda(x,\alpha)$
is $K$-reducible.
  We apply this  to
show that for any fixed $n\ge 10$ and for any number field $K$, all but
finitely many
$K$-specializations of the degree $n$ generalized Laguerre polynomial
$
L_n^{(t)}(x)
$
are $K$-irreducible and have Galois group $S_n$.  In conjunction with the
theory of complex multiplication, we also show that for any $K$ and for any
$
n\ge 53
$,
all but finitely many of the $K$-specializations of the modular equation
$
\Phi_n(x, t)
$
are $K$-irreducible and have Galois group containing $PSL_2(\Z/n)$.
\end{abstract}

\maketitle

\tableofcontents

\section{Introduction}

Let $K$ be a number field.  Consider a polynomial
$
\lambda(x,t)\in K[x,t]
$
which is non-constant in each of $x$ and $t$; it 
can be viewed as a one-parameter
family of  $K$-polynomials in $x$.  If
$
\lambda
$
is irreducible in $K[x, t]$, the Hilbert irreducibility theorem furnishes
infinitely many
$
\alpha\in K
$
for which 
$\mathcal \lambda(x, \alpha)$ is $K$-irreducible.   It is then natural
to study the set of $\alpha\in K$ with reducible specialization.  
These exceptional sets are thin sets \cite[$\S 9.6$]{serre}, and the example
$
x^n - t
$
shows that they can be infinite.  Using techniques from
diophantine analysis, Fried \cite{fried2} bounded the number of exceptional
specializations of bounded height.  
Exceptional sets for concrete families have also been examined;
for example the irreducibility and Galois group of the
Generalized
Laguerre polynomial
\begin{equation}
L_n^{(t)}(x)
=
\sum_{j=0}^n (-x)^j \binom{n}{j} \prod_{k=j+1}^n (t+k).
                \label{lag}
\end{equation}
for various rational values of the parameter $t$ were studied by
Schur (\cite{schur1}, \cite{schur}); more recently, Feit \cite{feit} used them
to solve the inverse Galois problem over $\Q$ for certain
double covers of the alternating group $A_n$.
See also 
\cite{gow}, \cite{hajir}, \cite{sell}, \cite{hajir2}, for other
related results.  Note that in the papers just cited,
the focus is primarily on a related,
but different, question from the one we began with, namely that of
irreducibility and Galois properties
of
$
L_n^{(\alpha_n)}(x)
$
for suitable sequences 
$
\{\alpha_n\}_n
$.  For example, the case $\alpha_n=-1-n$ corresponds to the truncated
exponential polynomial studied by Schur \cite{schur1}.
For the latter type of question, the $p$-adic Newton polygon is
a powerful tool.  For example, in Filaseta-Lam 
\cite{filaseta-lam} it is shown
that if we fix
$
\alpha\in \Q - \Z_{<0}
$,
then
$
L_n^{(\alpha)}(x)
$
is $\Q$-irreducible for $n$ sufficiently large, while in Filaseta-Trifonov
\cite{filaseta-trifonov}, Grosswald's conjecture, to the effect that
$
L_n^{(-1-2n)}(x)
$
(i.e.~the $n$-th degree Bessel polynomial)
is $\Q$-irreducible for every $n$, is proved.  The Newton Polygon 
approach, however,
does not appear to be  well-suited to the problem under consideration here,
namely that of studying exceptional specializations of $L_n^{(t)}(x)$
for $n$ {\em fixed}.

\vsp

In this paper we investigate the exceptional set of a given
$
\lambda(x, t)
$ from
the algebro-geometric and group-theoretic points of view.   First, note
that
$
\lambda(x, t)
$
defines a $1$-dimensional subvariety
$
X_\lambda\subset \mathbf P_K^2
$.
To say that the specialization of $\lambda$ at $t=\alpha$ has a
$K$-rational root is to say that the fiber above $\alpha$ of the
projection-to-$t$ map has a $K$-rational point.  Say $X_\lambda$ is in
fact absolutely irreducible; then, by Faltings, at most finitely many
$K$-specializations of $\lambda$ have a $K$-rational root if
$X_\lambda$ has genus $\ge 2$.  More generally, a result of M\"uller
\cite{muller} leads to an irreducibility criterion for specializations
in terms of the genus of intermediate subfields of $K'/K(t)$ where
$K'$ is the Galois closure of $\lambda(x, t)$ over the function field
$K(t)$ (cf.~also the related results of D\`ebes and Fried
\cite{fried}).  In sections \ref{sec:rational}--\ref{sec:spec}, we develop
and refine tools for applying this criterion.
In section \ref{sec:general}, we apply these
to study
$
L_n^{(t)}(x)
$.
The recursive properties of
$
L_n^{(t)}(x)
$
allow us to analyze the geometry
of the corresponding curve as well as the ramification behavior of the
projection-to-$t$ map.  By utilizing, in addition,
information about maximal subgroups
of the symmetric group $S_n$,  we obtain the following result.

\vsp

\begin{thm}
	\label{thm:irr}
Let $K$ be a number field.

{\rm (a)}  Fix $n\ge 5$.  Then for all but finitely many $\alpha\in K$,
$
L_n^{(\alpha)}(x)
$
is $K$-irreducible and its Galois group (over $K$) contains $A_n$.   For fixed
$n\ge 10$, this Galois group is
exactly $S_n$ except for finitely many $\alpha \in K$.

{\rm (b)}  Let $R$ be a finitely generated subring of $K$.  If $n\ge 6$, then
for all but finitely many $\alpha\in R$, the Galois group over $K$ of
$
L_n^{(\alpha)}(x)
$
is exactly $S_n$.
\end{thm}

\vsp

\begin{remark}
Note that Theorem \ref{thm:irr} is optimal in two ways.  First, for
$
6\le n\le 9
$,
the set of $\alpha\in K$ for which the discriminant of $L_n^{(\alpha)}(x)$ is a
square in $K$ turns out to be parameterized by a curve of geometric genus one,
so
for suitable $K$ there are infinitely many specializations with even Galois
group.  And when $n=5$, the square discriminants are parameterized by a curve
of geometric genus zero, so there are fields $K$ and finitely generated
subrings $R$
of $K$ over which there are infinitely many even specializations.
Second, 
$
L_4^{(t)}(x) = 0
$
is a model (cf. \cite{hajir2}) 
of the elliptic curve $384H2$ in Cremona's table.
This curve has
Mordell-Weil rank $1$ over $\Q$, so over any number field $K$ there are
infinitely many $\alpha\in K$ for which $L^{(\alpha)}_4(x)$
has a $K$-rational linear factor.  However, 
the exceptional set in Theorem \ref{thm:irr} is captured by rational
points on curves
of high geometric genus, so it would be difficult to make the Theorem
\textit{effective}.
\end{remark}

\vsp

Before we develop the tools necessary for proving Theorem ref{thm:irr},
we illustrate the use of M\"uller's criterion by applying it to 
another well-studied polynomial, namely the
modular polynomial $\Phi_n(x,j)$.
This monic $\Z$-polynomial plays a central role in the theory of elliptic
curves; it is  determined up to a scalar multiple by the property that
two elliptic
curves over $\C$ with $j$-invariants $j_1, j_2$ are related by a cyclic
$n$-isogeny if and only if 
$
\Phi_n(j_1, j_2) = 0
$.
It is 
irreducible over $\C(j)$, and its Galois group over $\Q(j)$ is
$
PGL_2(\Z/n)
$.

\vsp

For any integer $n>1$ and any prime $p$, define
$$
\Q_{p, n} =
\left\{
  \begin{array}{llllll}
    \text{unique quadratic extension of $\Q$ of conductor $p$}
    &
    \text{if $p>2$ and $p|n$,}
    \\
    \text{unique biquadratic extension of $\Q$ of conductor $8$}
    &
    \text{if $p=2$ and $8|n$,}
    \\
    \text{unique quadratic extension of $\Q$ of conductor $4$}
    &
    \text{if $p=2$ and $4||n$,}
    \\
    \Q                                  & \text{otherwise.}
  \end{array}
\right.
$$
For any number
field $K$ and any $n>1$, denote by $\tilde{K}_n$ the compositum of $K$ 
with all
$
\Q_{p, n}
$ 
as $p$ runs over the prime divisors of $n$;
note that this is a finite extension of $K$.

\vsp

\begin{thm}
	\label{thm:modular}
Let $n\ge 53$, and let $K$ be a number field.  Then for all but finitely many
$\alpha\in \tilde{K}_n$,
$
\Phi_n(x, \alpha)
$
is $K$-irreducible, and its Galois group over $\tilde{K}_n$ is $PSL_2(\Z/n)$.  If
$n$ is a prime then it suffices to take $n\ge 23$.
\end{thm}

\begin{remark}
Theorem \ref{thm:modular} is close to optimal in the $n$-aspect;
cf.~Remark \ref{rem:52}.  However, as in the discussion
following Theorem \ref{thm:irr}, it would be
difficult to make Theorem \ref{thm:modular} effective.
\end{remark}

\vsp

We will describe our strategy via M\"uller's criterion in section
\ref{sec:rational}, after we establish some notation.  To apply this
criterion to specializations of $\Phi_n$, in section \ref{sec:modular}
we investigate the algebraic closure of $\Q$ in the function field
defined by $\Phi_n$, and we study the genus of Riemann surfaces
defined by congruence subgroups.  In sections \ref{sec:rh} and
\ref{sec:spec}, we develop the technical tools needed for carrying out
the strategy outlined in section \ref{sec:rational}.  In section
\ref{sec:general}, we implement this plan for the Generalized Laguerre
Polynomial after first establishing several geometric properties of
the projective plane curve $\nl$ defined by $L_n^{(t)}(x)=0$.
Specifically, let $ \iota_n: \nl\rarr\mathbf P_K^1 $ be the branched
cover defined by the projection-to-$t$ map.  Then
\begin{itemize}
\item[(i)]
$K$ is algebraically closed in the splitting field of $L_n^{(t)}$ over
$K(t)$;
\item[(ii)]
the (geometric) Galois group of $\iota_n$ is $S_n$;
\item[(iii)]
$L_n^{(t)}(x)$, as a polynomial in $x$, has discriminant which is non-constant
in
$t$;
\item[(iv)]
$\nl$ has no affine singular points, and
\item[(v)]
$\iota_n$ has several ``simple'' branch points of index close to
$n$.
\end{itemize}
In (v), a {\em simple branch point of index $e$} is one
whose fiber consists of a number (possibly $0$) of multiplicity one
points together with a single ramified point (of multiplicity $e$).  
The cover defined by the degree $n$ Generalized Laguerre Polynomial
has one simple branch point of every index between $2$ and $n$: we use
the four of highest index, which suffices in our analysis for all $n\geq 6$.
As the calculations in section \ref{sec:general} will show, the proof of
Theorem \ref{thm:irr} extends readily to other one-parameter families
of polynomials satisfying
properties (i)-(v) (as long as their degree is large with respect to the
precise form taken by condition (v)).  
On the other hand, given an arbitrary $\lambda(x, t)$
which is irreducible over $K(t)$, in general we cannot expect all but
finitely many of its $K$-specializations to be $K$-irreducible, let
alone having the same Galois group as $\lambda(x,t)$ over $K(t)$ ---
the subvariety $X_\lambda$ mentioned just before the statement of
Theorem \ref{thm:irr} could, for example, have genus $\le 1$.  In
section \ref{sec:simple} we will analyze this situation further in the
case of ``simple branched covers,'' i.e. where all the branch points
are simple of index $2$.

\vsp

\section{Rational specializations}
      \label{sec:rational}

We first establish some notation and hypotheses which will be maintained
throughout.
Let $K$ be a field of characteristic $0$, finitely generated over $\Q$.  
Fix an algebraic closure $\ov{K}$
of $K$.  Denote by $K_0$ the function field $K(t)$.  Fix
$
\lambda(x, t)\in K[x,t]
$
so that $\lambda$ has degree $n>0$ in $x$ and is irreducible over $K_0$.  Then
$
K_1 := K[x]/(\lambda(x, t))
$
is a degree $n$ extension of $K_0$.  Let $K'/K_0$ be a Galois closure of
$K_1/K_0$, and write
$
G_\lambda = \gal(K'/K_0)
$.
By \cite[p.~123]{serre}, the Galois group of
$
\lambda(x, \alpha)
$
over $K$ is a subgroup of $G_\lambda$ for any $\alpha\in K$, and by
\cite[Prop.~9.2]{serre}, there are infinitely many $\beta_0\in K$ for which
this Galois group is exactly
$
G_\lambda
$.

\vsp  From now on, suppose that
\begin{itemize}
\item[(i)]
$K$ is algebraically closed in $K'/K_0$.
\end{itemize}
Then \cite[Remark II.2.5]{joe1} implies that every
intermediate subfield $E$ of $K'/K_0$ is the function field of a smooth
projective curve $X_E$ over $K$, and if
$
E\subset E'
$
are two such subfields, then there exists a $K$-morphism $X_{E'}\rarr X_E$
of degree
$
[E': E]
$.
We write $g(X_E)$ for the genus of $X_E$.
By Galois theory, intermediate fields $E$ of $K'/K_0$ 
are in bijective correspondence with subgroups
$\E=\gal(K'/E)$ of $G_\lambda$.  
To simplify the exposition, we abbreviate the phrase
`all but finitely many $\alpha\in K$' by 
$
\af.
$

\vsp

\begin{prop}
	\label{prop:muller0}
Let $K'/K_0$ be as above, and consider a polynomial
$f\in K[x, t]$ which is irreducible over $K_0$ but splits
completely into linear factors over $K'$.
Suppose for every intermediate subfield $E$ of $K'/K_0$ such that
$f$ is reducible over $E$,
we have $g(X_E)>1$.
Then $f(x, \alpha)$ is $K$-irreducible
for $\af$.
\end{prop}
\begin{proof} 
This is probably well-known to the expert; for a convenient reference see
M\"uller \cite[Prop.~4.20]{muller}.  A method of proof is also indicated
in \ref{subsec:fiber}.
\end{proof}

\vsp

For any $\alpha\in K$, the Galois group of
$
\lambda(x, \alpha)
$
over $K$ is a subgroup of $G_\lambda$, and we are interested in finding
conditions on
$\alpha$ under which  $\lambda(x,\alpha)$ is not only $K$-irreducible, 
but also has Galois group coinciding with the full $G_\lambda$.  Here is our
strategy: suppose the splitting field of some ``test-polynomial''
$
f(x, t)\in K[x,t]
$
is contained in $K'$; then the splitting field of 
$
f(x, \alpha)
$
over $K$ is contained in that of $\lambda(x, \alpha)$.  So if 
$
f(x, \alpha)
$
is $K$-irreducible, then the degree of the splitting field of
$
\lambda(x, \alpha)
$
over $K$ would be divisible by the degree of $f(x, \alpha)$.  By running
through an appropriate collection of $f$ (e.g. the polynomials
$\Lambda_j$ introduced in \ref{sec:spec}), we can then hope to show that
$
\#G_\lambda 
$
divides the degree of the splitting field of $\lambda(x, \alpha)$ over $K$,
whence the Galois group of
$
\lambda(x, \alpha)
$
over $K$ must be $G_\lambda$.   To study the irreducibility of the
specializations
$
f(x, \alpha)
$ 
we use Proposition \ref{prop:muller0}, which reduces the problem to estimating
the genus of $X_E$ as we run through intermediate subfields $E$ of $K'/K_0$.

\vsp

\section{Modular equations}
	\label{sec:modular}

By \cite[p.~55]{lang:cm}, the modular polynomial
$
\Phi_n(x, j)\in \Z[x, j]
$
is irreducible over $\C(j)$.
%
%
We now apply the strategy developed
in the last section to study specializations of $\Phi_n$.  Denote by 
$L_n$ the splitting field of $\Phi_n$ over $\Q(t)$.  Recall the definition of
$\Q_{p, n}$ and $\tilde{K}_n$ immediately preceding the statement of Theorem
\ref{thm:modular}.


\vsp

\begin{lem}
	\label{lem:cm}
The algebraic closure of $\Q$ in $L_n/\Q(t)$ is $\tilde{\Q}_n$.
\end{lem}

\begin{proof}
As a coarse moduli scheme, the open modular curve $Y_0(n)$ classifies 
isomorphism classes 
$
(E\rarr E')
$
of pairs of elliptic curves related via a cyclic $n$-isogeny.  Over the
complex numbers, such a
pair is completely determined by the $j$-invariants of $E$ and $E'$.  Thus
the
\textit{complex}
points of $Y_0(n)$ are canonically identified with the \textit{complex} points
of the affine plane curve defined by
$
\Phi_n(x, j)=0
$.
Under this identification, the projection-to-$j$ map from this complex plane
curve corresponds precisely to the branched cover
$
\pi_0(n): Y_0(n) \rarr Y_0(1)
$
coming from the inclusion
$
\Gamma_0(n)\subset SL_2(\Z)
$.
The smallest regular branched cover containing $\pi_0(n)$ is then the cover
$
\pi(n): Y(n)\rarr Y(1) = Y_0(1)
$
corresponding to the inclusion
$
\Gamma(n) \subset SL_2(\Z)
$.
In particular, the deck transformation group of $\pi(n)$ is
$$
PSL_2(\Z)/ ( \Gamma(n)/\pm I) \simeq PSL_2(\Z/n).
$$
It follows that the \textit{geometric} Galois group of $\Phi_n$ is
$PSL_2(\Z/n)$.  But Macbeath \cite{macbeath} showed that
$
\gal(L_n/\Q(t)) \simeq PGL_2(\Z/n)
$,
so the algebraic closure of $\Q$ in $L_n/\Q(t)$ is the compositum of $\Q(t)$
with a Galois extension $L(n)/\Q$ with Galois group
\begin{eqnarray}
PGL_2(\Z/n)/ PSL_2(\Z/n)
&\simeq&
\prod_{p|n} PGL_2(\Z/p^{e_p})/ PSL_2(\Z/p^{e_p})
         \hspace{20pt}
	 \text{where $p^{e_p}|| n$}             	 \nonumber
\\
&\simeq&
\prod_{\substack{p|n\\p>2}} (\Z/2)
\times
\biggl\{
  \begin{array}{llllll}
    \Z/2\times \Z/2    &    \text{if $8|n$}
    \\
    \Z/2               &    \text{if $4||n$}
    \\
    \{ 1 \}            &    \text{otherwise}
  \end{array}
\biggr\}.
            \label{gal}
\end{eqnarray}
If $m|n$ then $L_m\subset L_n$, hence $L(m)\subset L(n)$, so to prove the
Lemma we are reduced to showing that for any prime power
$p^e>1$,
\begin{equation}
L(p^e) = \Q_{p, p^e}.         \label{pe}
\end{equation}

\vsp

For any $\alpha\in\Q$ and any $n>1$, the splitting field of
$
\Phi_n(x, \alpha)
$
over $\Q$ also contains $L(n)$.  Take $\alpha\in\Q$ to be one of the thirteen
$j$-invariants over $\Q$ corresponding to CM elliptic curves over $\Q$,
say
$
\alpha = j(\tau)
$.
Denote by $k_\alpha/\Q$ the corresponding complex quadratic field.  By the
`First Main Theorem' of complex multiplication \cite[Thm.~11.1]{cox},
$
k_\alpha( j(n\tau) ) 
$
is the ring class field of $k_\alpha$ of conductor $n$, hence
$
L(n)\subset k_\alpha( j(n\tau) ) 
$.
In particular, $L(n)/\Q$ is unramified outside of the prime divisors of $n$ and
of the discriminant of
$
k_\alpha/\Q
$.
If
$
j(\tau') = \alpha'\in\Q
$
is another CM $j$-invariant over $\Q$, then 
$
L(n)\subset k_\alpha( j(n\tau) )  \cap k_\alpha'( j(n\tau') )
$.
We may choose $\alpha'$ so that $k_\alpha$ and $k_\alpha'$ have coprime
discriminants, whereby $L(p^e)/\Q$ is unramified outside $p$.  On the
other hand,
(\ref{gal}) says that
$
L(p^e)/\Q
$
is quadratic if $p>2$ or $p^e=4$, and that it is biquadratic if $8|p^e$.
Recalling the definition of
$
\Q_{p, n},
$
we get (\ref{pe}) except when $p^e=4$.
To treat this remaining case we actually need to determine these ring class
fields.

\vsp

Set
$
\omega = \frac{1+\sqrt{-7}}{2}
$,
and take $\alpha = j(\omega)\in \Q$, so $k_\alpha=\Q(\omega)$.  
The conductor of the extension
$
k_\alpha( \sqrt{-1}) /k_\alpha
$
clearly divides $4\Z[\omega]$.  On the other hand, by \cite[Thm.~7.24]{cox}
the ring class field of $k_\alpha$ of conductor $4\Z[\omega]$ is a quadratic
extension of $k_\alpha$, so this ring class field is precisely
$
k_\alpha(\sqrt{-1})
$.
Recalling (\ref{gal}), we see that $L(4)/\Q$ is a quadratic extension in
$
\Q(\omega, \sqrt{-1})
$
unramified outside $2$, and (\ref{pe}) follows for $p^e=4$.
\end{proof}

\vsp

Rademacher conjectured that there are only finitely many congruence subgroups
with corresponding modular curve of genus zero (cf.~\cite{knopp}).
Dennin \cite{dennin:n} proved the stronger result that for any integer $g$,
there are at most
finitely many $n$ for which $PSL_2(\Z/n)$ contains a subgroup of genus $\le g$.
Cummins and Pauli \cite{cummins} recently tabulated all such subgroups for
$g\le 24$, from which we deduce the following result.

\vsp

\begin{lem}[Cummins-Pauli]
	\label{lem:genus}
If $n\ge 53$, then every proper subgroup of $PSL_2(\Z/n)$ has
genus $\ge 2$.  If $n$ is a prime, the same conclusion 
holds for  $n\ge 23$.
       \qed
\end{lem}

\vsp

\begin{proof}[Proof of Theorem \ref{thm:modular}]
Thanks to Lemma \ref{lem:cm}, the discussion in section
\ref{sec:rational} is applicable to
$
\Phi_n
$
over $\tilde{K}_n$ for any number field $K$.

\vsp

Let $\pi_n$ be a primitive element for the extension
$
\tilde{K}_n L_n(t)/\tilde{K}_n(t)
$,
and let $f_n(x, t)$ be the minimal polynomial of $\pi_n$ over $\tilde{K}_n(t)$.
Then $f_n$ is irreducible over 
$
\tilde{K}_n(t)
$,
by construction.  So if $n$ is as in Lemma \ref{lem:genus}, then
Proposition \ref{prop:muller0} and this Lemma together imply that for
$
\alpha\in_{\text{af}} \tilde{K}_n
$,
the specializations of $f_n$ and of $\Phi_n$ at $t=\alpha$ are both
$\tilde{K}_n$-irreducible.  If we write 
$
F_n(\alpha)
$
for the splitting field of $\Phi_n(x, \alpha)$ over $\tilde{K}_n$, 
then that means
$
[F_n(\alpha): \tilde{K}_n]
$
is divisible by
$
\deg f_n = [\tilde{K}_n L_n(t): \tilde{K}_n(t)] = \#PSL_2(\Z/n)
$,
and Theorem  \ref{thm:modular} follows.
\end{proof}

\vsp

\begin{remark}
  \label{rem:52}
The fact that \textit{every} non-trivial intermediate subfield of 
$
\tilde{K}_n L_n(t)/\tilde{K}_n(t)
$
has genus $\ge 2$ for $n\ge 53$ significantly simplifies our search for the
`test polynomial' $f$ in Proposition \ref{prop:muller0}. The modular curve
$X_0(n)$
has genus $\le 1$ for $n\le 21$ and for
$
n \in \{24, 25, 27, 32, 36, 49\},
$
so by the discussion
immediately preceding Theorem \ref{thm:irr}, for these $n$ the modular
equation has infinitely many reducible specializations over suitable $K$.
To analyze the remaining values of $n\le 52$ we could search for test
polynomials $f$ which remain
irreducible over intermediate subfields of genus $\le 1$.  We will not pursue
this issue here, but in section \ref{sec:spec} we will study the same
problem
for specializations of $S_n$-extensions by using a family
of $[(n-1)/2]$ test polynomials
$
\Lambda_j(x, t)
$.
\end{remark}

\vsp

\section{A Riemann-Hurwitz estimate}
	\label{sec:rh}

We now return to the general setup in section \ref{sec:rational}.  To apply
Proposition \ref{prop:muller0}, 
we need to be able to estimate the genus of certain intermediate
subfields of $K'/K_0$.
To do that we will apply the Riemann-Hurwitz formula to the
cover
$
\xi_E: X_E\rarr \po
$
corresponding to the field inclusion $K_0 \subset E$.  Since we do not have
any explicit model for $X_E$, we will take an algebraic approach.
Thanks to hypothesis (i) in section \ref{sec:rational}, in order to determine
the ramification of the \textit{geometric} cover
$
X'\rarr\po
$
it suffices to determine the \textit{algebraic} ramification behavior of
integral extensions of Dedekind domains corresponding to this geometric cover.

%
%
%
%
%
%

\vsp

Denote by $ \be\subset\po $ the branch locus of the
projection-to-$t$ map for $\lambda$.  Then $\xi_E$ is unramified
outside $\be$.  Fix affine open sets on $X_E$ and $X'$ which contain
every fiber of $\xi_E$ and $X'\rarr \po$ above $\be$, and denote by
$\oo_E$ and $\oo'$ their respective affine coordinate rings.  Write
$\oo_0$ for the affine coordinate ring of the affine line in $\po$.
Let $\m_\nu$ (or just $\m$ if $\nu$ is fixed) be the maximal ideal in
$\oo_0$ corresponding to a given $ \nu\in\be $.  We let $e_\nu=e(\M/\m)$
be the ramification index of $\M$ in the Galois cover $K'/K_0$, where
$\M$ is an arbitrary prime of $\oo'$ dividing $\m\oo'$.

%

\begin{defn}\label{def1}
\begin{itemize}
\item[(a)] For a positive integer $\delta$ and a branch point $\nu\in \be$ 
corresponding to an ideal $\m \in \oo_0$, let 
$$
c_\delta(\nu)= c_\delta(\m)= \sum_
{\substack{\n|\m\oo_E\\e(\n/\m) = \delta}}
f(\n/\m),
$$ 
be the sum of the residual degrees of distinct $\oo_E$-primes
$\n$ of ramification index $\delta$ over $\m$.

\item[(b)] For $\nu \in \be$ corresponding to an ideal $\m \in \oo_0$, let
$$\Delta(\nu) = \Delta(\m) = \sum_
{{\n|\m\oo_E}}
(e(\n/\m)-1)f(\n/\m)
$$
be the $\nu$-component of the discriminant of $E/K_0$.

\item[(c)] For an integer $e>1$, let $d(e)$ be the least prime divisor of $e$.
\end{itemize}
\end{defn}



\begin{lem}\label{r-h}
With the notation and hypotheses as in section \ref{sec:rational},
if $E$ is an intermediate field of $K'/K_0$ corresponding to a
subgroup $\E=\gal(K'/E)$ of $G_\lambda=\gal(K'/K_0)$, and $V$ is any subset
of $\be$, then 
\begin{equation}
               \label{cd1}
g(X_E) \geq 
1 + \frac{[G: \E]}{2}
        \Bigl(
                -2
                +
                \sum_{\nu\in V} \Bigl( 1- \frac{1}{d(e_\nu)} \Bigr)
        \Bigr)
        -
        \frac{1}{2}
        \sum_{\nu\in V} c_1(\nu) \Bigl( 1 - \frac{1}{d(e_\nu)} \Bigr).
 \end{equation}
\end{lem}
\begin{proof}
First, note that 
\begin{equation}
\sum_{1\leq\delta|e_\nu} c_\delta(\nu) \delta = [E: K_0] = [G_\lambda: \E].
	\label{efg}
\end{equation}
For each $\nu \in \be$, we have from Definition \ref{def1}, 
\begin{eqnarray}
\Delta(\nu)
&=&
\sum_{1\leq \delta | e_\nu}  c_\delta(\nu)(\delta - 1)
	\nonumber
\\
&=&
\sum_{1<\delta| e_\nu}
       c_\delta(\nu)
       \Bigl( 1 - \frac{1}{\delta} \Bigr) \delta
	\nonumber
\\
&\ge&
\Bigl( 1 - \frac{1}{d(e_\nu)} \Bigr)
       \sum_{1<\delta|e_\nu} c_\delta(\nu) \delta
	\nonumber
\\
&\ge&
\Bigl( 1 - \frac{1}{d(e_\nu)} \Bigr)
\sum_{1\leq \delta | e_\nu}
       c_\delta(\nu) \delta
-
\Bigl( 1 - \frac{1}{d(e_\nu)} \Bigr) c_1(\nu) 
	\nonumber
\\
&\geq&
[G_\lambda: \E] \Bigl( 1 - \frac{1}{d(e_\nu)} \Bigr)
-
c_1(\nu) \Bigl( 1 - \frac{1}{d(e_\nu)} \Bigr)
        \hspace{20pt}	\text{by (\ref{efg}).}
        \label{more}
\end{eqnarray}

By Riemann-Hurwitz for $E/K_0$, \cite[Theorem 7.16]{rosen}, we have 
$$
g(X_E) - 1 =  [E:K_0](0-1) + \frac{1}{2} \sum_{\nu \in \be} \Delta(\nu).
$$
Since $\Delta(\nu)>0$, we have, for any subset $V\subseteq B_\lambda$, 
\begin{eqnarray}
g(X_E)
&\ge&
1 - [G_\lambda: \E] + \frac{1}{2}\sum_{\nu\in V}
		\Delta_\nu
	\nonumber
\\
&\ge&
1 + \frac{[G_\lambda: \E]}{2}
	\Bigl(
	        -2
		+
		\sum_{\nu\in V} \Bigl( 1- \frac{1}{d(e_\nu)} \Bigr)
	\Bigr)
	-
	\frac{1}{2}
	\sum_{\nu\in V} c_1(\nu) \Bigl( 1 - \frac{1}{d(e_\nu)} \Bigr)
\hspace{10pt}
\text{by (\ref{more}).}
\nonumber
\end{eqnarray}
\end{proof}
\begin{remark}
Note that the bound (\ref{cd1}) is useful only when
$c_1(\nu)$ is fairly small for all $\nu \in V$, so in using (\ref{cd1}),
it's often useful to take $V$ to be a proper subset of $\be$. 
Moreover, the inequality (\ref{cd1}) is in fact strict if $V$ is a
{\em proper} subset of $\be$ since $\Delta(\nu)>0$ for $\nu\in V$.
\end{remark}

In view of Proposition \ref{prop:muller0}, our task will be to show that the
right hand side of (\ref{cd1})
is $>1$ when a given
$
f(x, t)\in K[x, t]
$
is reducible over $E$.  
For our application
to Generalized Laguerre Polynomials, this will be easy to arrange by
taking $V$ to be an appropriately small subset of $\be$.

\vsp

We now turn to the task of bounding $c_1=c_1(\nu)$ from above,
where, for the remainder of this section,
$\nu \in \be$ is a fixed branch point, with corresponding ideal
$\m=\m_\nu$ of $\oo_0$.
Fix also
a prime $\M \subset \oo'$ lying over $\m$,
with corresponding decomposition group $D=\{ \sigma\in G :
\M^\sigma = \M\},$ and inertia group $I=I(\M/\m)$.  Let $T$
be a subset of $G=G_\lambda=\gal(K'/K_0)$ such that 
\begin{equation}\label{dcd}
G=\coprod_{\tau\in T} \E \tau D
\end{equation}
is the decomposition of $G$ into
disjoint double cosets, where 
$\E=\gal(K'/E)$ is the subgroup fixing $E$.

As is clear from Lemma \ref{r-h}, it will be important to keep
track of the primes $\n$ of $\oo_E$ dividing $\m$ and especially their
ramification indices $e(\n/\m)$.  That these can be described nicely
in terms of the double coset decomposition (\ref{dcd}) is a useful
fact (we learned from Tate) for which we were not able to find a
suitable reference, so we give the details.  For each $\sigma \in G$,
let $\n_\sigma$ be the prime $\M^\sigma \cap \oo_E$ of $\oo_E$ lying
under $\M^\sigma$.  Let $I_\sigma \subseteq D_\sigma$ be the inertia
and decomposition groups of $\M^\sigma/\m$, respectively.  They
satisfy $D_\sigma=\sigma D \sigma^{-1}$ and $I_\sigma=\sigma I
\sigma^{-1}$.  In the extension $K'/E$, the inertia and decomposition
groups for $\M^\sigma/\n_\sigma$ are simply $I_\sigma \cap \E$ and
$D_\sigma \cap \E$, respectively.  For the ramification indices of
$\M/\m, \M/\n_\sigma$, and $\n_\sigma/\m$, let us put
$$
e = e(\M/\m), \qquad e_\sigma'=e(\M/\n_\sigma), 
\qquad e_\sigma=e(\n_\sigma/\m),
$$
and similarly for the residual degrees, we put
$$
f = f(\M/\m), 
\qquad f_\sigma'=f(\M/\n_\sigma), \qquad f_\sigma=f(\n_\sigma/\m).
$$
By multiplicativity in towers for these invariants, we have 
\begin{equation}\label{ef}
e_\sigma e_\sigma'=e, \qquad f_\sigma f_\sigma' = f.
\end{equation}

\begin{lem}\label{tate}
With the notation introduced above,
\begin{enumerate}\item[(a)] The distinct primes of $\oo_E$ dividing $\m$ are
those induced by $\M^\tau$ for $\tau \in T$.  In other words, 
we have $\n_\sigma=\n_{\sigma'}$ if and only if 
$\E \sigma D = \E \sigma' D$.

\item[(b)]  For $\sigma \in G$, 
we have $$e_\sigma f_\sigma = [ \sigma D \sigma^{-1} : \E \cap
\sigma D \sigma^{-1}], \qquad e_\sigma = [ \sigma I \sigma^{-1} : \E \cap
\sigma I \sigma^{-1}].$$

\end{enumerate}

\end{lem}
\begin{proof}
Let $w$ be the valuation of $\oo'$ corresponding to $\M$.  For $\alpha
\in \oo'$, we have
$|\alpha|_{\sigma w}= |\sigma^{-1}\alpha|_w$.  If $\E
\sigma D = \E \sigma' D$, we can write $\sigma' = h \sigma g$, with
$h\in \E, g \in D$.  For $\alpha \in \oo_E$, we compute
$$ |\alpha|_{\sigma' w} = |\alpha|_{h\sigma gw}=|\alpha|_{h\sigma w}
=|h^{-1} \alpha|_{\sigma w} = |\alpha|_{\sigma w}.$$ Thus, $\sigma w$
and $\sigma' w$ induce the same valuation on $\oo_E$, i.e. $\n_\sigma =
\n_{\sigma'}$.  Conversely, suppose $\n_\sigma = \n_{\sigma'}$,
i.e. the set of primes of $\oo'$ lying over $\n_\sigma$ includes $\M^{\sigma'}$ 
as well as $\M^\sigma$.  Since $\E=\gal(K'/E)$ acts transitively on this
set, there exists $h \in
\E$ such that $\M^{h \sigma'}=\M^{\sigma}$, i.e. $\sigma^{-1}h \sigma'
\in D$. Therefore, $\E\sigma D=\E \sigma' D$.  This proves (a).  
 We have $ef=\#I_\sigma f= \#D_\sigma$ and
$e_\sigma'f_\sigma'=\#(I_\sigma\cap \E)f_\sigma'=\#(D_\sigma\cap \E)$,
so we get (b) by multiplicativity in towers (\ref{ef}).
\end{proof}

Define
$$
Y = \{ \sigma \in G:  \sigma I \sigma^{-1} \subset \E \}.
$$
For the application to Riemann-Hurwitz,
we'll need to estimate $c_1$.  We proceed as follows.

\begin{lem}
  \label{lem:far}
If $a \in Y$, then $\{ b \in Y: \E aI = \E bI \} = \E a$.  We have 
$ c_1 = \#Y/\# \E$.
\end{lem}

\begin{proof}
We first make a remark that simplifies the calculation.  Note that if
we compose our fields $K_0\subset E\subset K'$ with a finite extension
$\tilde{K}$ of the constant field $K$ that splits $\M$, then
$c_\delta(\m)$ remains unchanged, since each prime $\n_\sigma$ of $E$
of residual degree $f_\sigma$ splits in $E\tilde{K}$ into $f_\sigma$
primes of residual degree $1$ with the same inertia group
$I_\sigma\cap \E$.  In fact, the genus calculation we are performing
is a purely geometric one, so we could have simply assumed from the
outset that the constant field $K$ is algebraically closed.

Either way, we take $\M/\m$ as above and assume
without loss of generality, that $f(\M/\m)=1$, i.e. $I=D$.

By Lemma \ref{tate}, for any $\sigma\in G$,
$
e(\n_\sigma/\m)=1
$
if and only if $\sigma I\sigma^{-1}\subset \E$.  Thus
\begin{equation}
c_1 = \#\{ \E \sigma I:   \sigma I\sigma^{-1} \subset \E\}.        \label{c1}
\end{equation}

Note that
$
\E aI = \E bI
$
if and only if $b\in\E aI$.  Suppose $b\in Y$ and $b\in\E aI$.  Then
$
ba^{-1}\in\E aIa^{-1} \subset \E
$,
hence $b\in\E a$.  Conversely, suppose $b = h a$ with $h\in \E$.  Then
$$
b I b^{-1} = h a I a^{-1} h^{-1} \subset h \E h^{-1} = \E
$$
so $b\in Y$.
Finally, clearly $\E a \subset \E aI$ so $b \in \E a$ implies $b \in \E aI$.
Therefore, 
$Y$ is
a union of {\rm(}\!right{\rm)} cosets of $\E$, and the number of distinct
double
cosets $\E aI$ with $a \in Y$ is exactly $\#Y/\#\E$.
This completes the proof by (\ref{c1}).
\end{proof}

\vsp

Since we are working with function fields of characteristic $0$, all
ramification is tame, so the
inertia group $I$ is cyclic. We now specialize to the case where
$G=S_n$, and $I$ is generated by a {\em cycle} (under its natural action
on the roots of $\lambda$).  Of course, if $\#I$ is greater than
$n/2$, the latter condition holds automatically.
\begin{lem}\label{lem:c1bound}
If $\gal(K'/K_0)=S_n$ and $I$ is generated by an $m$-cycle, then 
\begin{eqnarray}
c_1
&=&
\frac{(\text{number of $m$-cycles in $\E$})}{\#\E}
\times
m(n-m)!
       \label{c11}
\\
&<&
m(n-m)!.
       \label{c12}
\end{eqnarray}
\end{lem}
\begin{proof}
Just as in the proof of the preceding Lemma, we may assume
that $I=D$.  Let
$
J = \{ sIs^{-1}\subset\E:  s\in G \}
$ be the set of subgroups of $\E$ which are $G$-conjugate to $I$.
Then
\begin{eqnarray}\label{yy}
\#Y
&=&
\sum_{I'\in J} \#\{ s\in G: sI's^{-1} = I' \}.
	\label{y1}
\end{eqnarray}
Any two $m$-cycles in $S_n$ are $S_n$-conjugate,
so
\begin{eqnarray}
\#J
&=&
\text{number of cyclic subgroups of $\E$ generated by an $m$-cycle}
\\
&=&
\text{(number of $m$-cycles in $\E$)}/\varphi(m).    \label{jj}
\end{eqnarray}
There are
$
n!/( m(n-m)! )
$
$m$-cycles in $S_n$, so for any $S_n$-conjugate $I'\subset\E$ of $I$,
\begin{eqnarray}
\#\{ s\in S_n: sI's^{-1} = I' \}
=
\frac{n!}{\# \text{orbit}_{S_n}(I')}
=
\frac{n!}{ \frac{n!}{m(n-m)!}/\varphi(m)}
=
m \varphi(m) (n-m)!.
                 \label{are}
\end{eqnarray}
The proof is complete once we combine 
(\ref{yy})-(\ref{are}) with Lemma \ref{lem:far}.
\end{proof}


We end this section with an elementary criterion which guarantees the
hypothesis of Lemma \ref{lem:c1bound} (on inertia being generated by a
cycle) to hold; the criterion will be easily verified for the
Generalized Laguerre Polynomial at all its branch points.

\vsp

Recall that $K_1/K_0$ is a root field for $\lambda$, i.e. $K_1\simeq
K_0[x]/(\lambda)$.

\begin{defn}
Let $\nu\in\be$ be a branch point of $\lambda$, with corresponding
maximal ideal $\m \subset \oo_0$.  Let $e>1$ be an integer.
We say that $\nu$ (or $\m$) is {\em simple of index $e$}
for $\lambda$ if
\begin{equation}
\m\oo_{K_1} = \n_0^{e} \n_1 \cdots \n_s,		\label{m}
\end{equation}
where $\n_0, \ldots, \n_s$ are pairwise distinct primes of $\oo_{K_1}$;
in other words, in $\oo_{K_1}$, there is a unique prime dividing
$\m\oo_{K_1}$ with non-trivial ramification index (equal to $e$).
\end{defn}

\begin{lem}
       \label{lem:ram}
Suppose $G=\gal(K'/K_0)=S_n$.  Let
$
\m\subset\oo_0
$
be a maximal ideal corresponding to a branch 
point  $\nu\in B_\lambda$, which is simple of index $e>1$.
Then,
for any
$
\M\subset\oo'
$
lying above $\m$, the inertia group
$
I=I(\M/\m)
$
has order $e$ and is generated by a cycle of
length $e$.
\end{lem}

\begin{proof}
Let $\E=\gal(K'/K_1)$.  The index $n$ subgroups in $S_n$ are
stabilizers of any one of the $n$ letters.  By reordering the roots if
needed, we can identify $\E \simeq S_{n-1}$ with the stabilizer of the
letter $n$.  Every element in $S_n$ is a product of disjoint,
non-trivial cycles.  This decomposition is unique once a labelling is
fixed, and two elements in $S_n$ are conjugate if and only if they
decompose into the same number of cycles of each length.  

Returning to the proof of the Lemma, suppose $\M$ is a prime of $\oo'$
whose restriction $\M\cap \oo_{K_1}$ is the unique prime $\n$ of
$\oo_{K_1}$ of ramification index $e>1$ over $\m$.  Let $I=I(\M/\m)$.
We may assume, as in the preceding lemmas, that composing with a
suitable finite extension of $K$, $\M/\m$ has degree $1$, i.e. $I=D$
(this disturbs neither the identification 
$G\simeq S_n$ nor the embedding  $I \hookrightarrow G$).  

Let $\gamma $ be a generator of the cyclic group $I$.  Write $ \gamma
= \gamma_1 \cdots \gamma_r $ for its decomposition into disjoint,
possibly trivial, cycles.  Since the $\gamma_i$ pairwise commute, we
may assume that the letter $n$ occurs in the cycle $\gamma_1$.  For
$1\leq i \leq r$, let $ a_i = \text{ord}(\gamma_i)\geq 1 $, and let
$a=\min \{ m \geq 1 : \gamma^m \in \E \}$.  On the one hand,
$\gamma^a$ generates $I\cap \E$, and, on the other hand, we have
$a=a_1$ (recalling our convention that $\E$ is the stabilizer of the
letter $n$).  By Lemma \ref{tate}, $e(\n/\m)=\#[I/(I\cap
\E)]=\#[\langle \gamma \rangle / \langle \gamma^a \rangle] = a $, thus
$\gamma_1$ has order $a=a_1=e>1$, since we took $\n$ to be the unique
prime of ramification index $e>1$ over $\m$.

It remains to show that the cycles $\gamma_2, \ldots, \gamma_r$ are
trivial, i.e.  $a_i=1$ for $i>1$.  We proceed by contradiction.  If
$a_2>1$, say, then, there exists $\sigma \in G$ such that $\sigma
\gamma_2 \sigma^{-1}$ is a cycle acting non-trivially on the letter
$n$.  Then, as before, $e(\n_\sigma/\m)=a_2>1$, so we get $a_2=e$ and
$\n_\sigma=\n$ by the assumption on the simplicity of the
ramification.  By Lemma \ref{tate}, therefore, $\sigma \in \E I$, say
$\sigma = \eta \theta$ with $\eta\in \E$ and $\theta \in I$.  Letting $x'=
\sigma x \sigma^{-1}$ for $x\in G$, we have $\gamma'= \gamma_1'
\gamma_2' \cdots \gamma_r'$ is the decomposition of $\gamma'$ into
disjoint cycles since conjugation preserves cycle structure.  But we
claim that $\gamma_1'$ and $\gamma_2'$ are not disjoint, as they both
act non-trivially on the letter $n$.  To see this, note that $\theta =
\gamma^b$ for some integer $b$, so $\theta \gamma_i \theta^{-1}=\gamma_i$
for $i=1,\ldots,r$.  On the other hand, since $\eta\in\E$, it fixes $n$, so
$\gamma_1'= \eta\gamma_1 \eta^{-1}$ and $\gamma_2'=\eta\gamma_2 \eta^{-1}$ are
both $e$-cycles that act non-trivially on $n$, hence are not disjoint.
This contradiction shows that $\gamma_2, \ldots,\gamma_r$ are all
trivial, so $I=\langle \gamma_1 \rangle$ is generated by an $e$-cycle,
hence has order $e$.
\end{proof}

\vsp

\vsp

\section{Specializations of $S_n$-covers}
                  \label{sec:spec}

In this section, we develop a strategy for applying Proposition
\ref{prop:muller0} to a geometric $S_n$-cover.
Namely, starting with an $S_n$-extension of function fields $K'/K_0$
as in section \ref{sec:rational}, in subsection \ref{subsec:inter} we
construct a family of polynomials $ \Lambda_j(x, t) \in K[x, t] $ with
splitting field contained in $K'$ (to which we will later apply
Proposition \ref{prop:muller0}).  In \ref{subsec:fiber}, we will give
a geometric interpretation in terms of fiber products for the curves
corresponding to these $\Lambda_j$ which we need for controlling the
genus of subfields of $K'$ cut out by a subgroup contained in $A_n$.
A reader who is interested in a proof of Theorem \ref{thm:irr} for 
$n\geq 10$ only, can skip \ref{subsec:fiber} entirely, as it will enter
the proof only for $6\leq n \leq 9$.

\vsp

\subsection{Distinguished subfields in $S_n$-extensions}
        \label{subsec:inter}

\mbox{ }

Let $\lambda(x, t)$ and $K'/K_0$ be as in section \ref{sec:rational};
in particular, recall the regularity hypothesis (i) introduced there.
Suppose further that
\begin{itemize}
\item[(ii)]
$G_\lambda\simeq S_n$, and
\item[(iii)]
$\lambda$, as a polynomial in $x$, has discriminant which is non-constant in
$t$.
\end{itemize}

\vsp
These two conditions actually recover the regularity of the cover,
at least when $n$ is not too small.
\begin{lem}
	\label{lem:alg}
Suppose $n\ge 5$.  Then

{\rm (a)}
$K$ is algebraically closed in $K'/K_0$, and

{\rm (b)}
$K'/K_0$ has a unique Galois subfield.  This subfield is quadratic over $K_0$.
\end{lem}

\begin{proof}
Fix an algebraic closure $\ov{K}$ of $K$.  Then $\ov{K}\cap K'$ is a Galois
subfield of the
$S_n$-extension $K'/K_0$.  Since $n\ge 5$,  the only non-trivial Galois
subfield in $K'/K_0$ is the unique quadratic subfield generated by the
square-root of the discriminant (with respect to $x$) of
$
\lambda(x, t)
$.
Invoke the discriminant condition on $\lambda$ and we are done.
\end{proof}

\vsp

The following result is standard.

\vsp

\begin{lem}
	\label{lem:plane2}
Let $X/K$ be a smooth projective curve,
and let $\xi: X\rarr\po$ be a non-constant $K$-morphism.  Then $X$ is 
$K$-birational to a plane curve $G(x, t)=0$ such that $\xi$ is the
projection-to-$t$ map.	\qed
\end{lem}

We now describe a distinguished collection of subfields in $K'/K_0$.
Fix a labelling of the roots of $\lambda(x, t)$ over $K_0$, giving an
identification of
$G_\lambda$ with the symmetric group $S_n$.  For
$
1\le j<n
$,
write $S_{n, j}$ for the subgroup
$
S_j\times S_{n-j} \subset S_n
$,
where $S_j$
permutes the first $j$ roots, and $S_{n-j}$, the remaining $n-j$ roots.
Denote by
\begin{itemize}
\item
$K_j$ the subfield of $K'/K_0$ fixed by $S_{n, j}$,
\item
$X_j$ the associated smooth projective curve over $K$, and
\item
$
\tilde{\phi}_{n,j}:   X_j\rarr \po
$
the $K$-branched cover corresponding to the extension $K_j/K_0$.
\end{itemize}
Lemma \ref{lem:plane2} furnishes a $K$-birational map taking $X_j$ to a
plane curve
$
\Lambda_j(x, t)=0
$
which is smooth above $t=\beta_0$, and such that 
$
\tilde{\phi}_{n, j}
$
is the projection-to-$t$ map.  Clearly we can take
$
\Lambda_1= \lambda
$ and do so.
Since $X_j$ is smooth, it is absolutely irreducible, hence so is
$
\Lambda_j(x, t)
$.
Thus we can apply Proposition \ref{prop:muller0} to $\Lambda_j$.

\vsp

\begin{lem}
	\label{lem:muller2}
Fix positive integers $n,j$ satisfying $n\geq 5$ and $j\in [1,n/2]$.
Suppose for every intermediate subfield $E$ of $K'/K_0$ over which
$
\Lambda_j(x, t)
$
is reducible, we have
$g(X_E)>1$.  
%
%
%
Then for $\af$, the
specialization
$
\lambda(x, \alpha)
$
is $K$-irreducible, and its splitting field has degree divisible by
$
\binom{n}{j}.
$
\end{lem}

\begin{proof}
As
$
\deg \tilde{\phi}_{n, j}
=
[K_j: K_0]
=
\#S_n/\#S_{n, j} = \binom{n}{j} \ge n
$,
and $n\ge 5$, Lemma \ref{lem:alg}(b) says that $K'/K_0$ is the
Galois closure of $K_j$; equivalently, $K'/K_0$ is the splitting field
of $\Lambda_j(x, t)$ over $K_0$.  But $K'$ is the splitting field of 
$
\lambda(x, t) = \Lambda_1(x, t)
$
over $K_0$, so by Proposition \ref{prop:muller0}, for $\af$ the splitting field
of
$
\lambda(x, \alpha)
$
contains the roots of $\Lambda_j(x, \alpha)$, and we are done.
\end{proof}

\vsp

For the proof of Theorem \ref{thm:irr}, we will employ the following
application of Proposition \ref{prop:muller0}.

\begin{thm}
	\label{thm:special}
Suppose $n\geq 7$ and
$
\Lambda_j(x, t)
$
satisfies the hypothesis in Lemma \ref{lem:muller2} for
each integer $j\in [1,n/2]$.  Then for
$\af$, the specialization
$
\lambda(x, \alpha)
$
is $K$-irreducible and has Galois group containing $A_n$.
\end{thm}

\begin{proof}
First, recall that $\Lambda_1 = \lambda$.  By Lemma \ref{lem:muller2},
$
\lambda(x, \alpha)
$
is $K$-irreducible for $\af$, hence its Galois group is a transitive subgroup
of $S_n$.  If $n\ge 8$, then there exists a prime $q$ with $n/2<q<n-2$
\cite[p.~370]{schur2}.  Necessarily $q$ divides
$
\binom{n}{k}
$
for some $1<k<n/2$, so by Lemma \ref{lem:muller2}, for $\af$ the specialization
$
\lambda(x, \alpha)
$
is $K$-irreducible, and $q$ divides the degree of its splitting field over $K$.
That means the Galois group of such a
$
\lambda(x, \alpha)
$
is a transitive subgroup of $S_n$ and has order divisible by $q$; a
theorem of Jordan \cite[Thm~5.6.2 and 5.7.2]{hall} then implies that this
Galois group contains $A_n$.

\vsp

For $n=7$,  Lemma \ref{lem:muller2} implies that for
$\af$, the Galois group of 
$
\lambda(x, \alpha)
$
is a transitive subgroup of $S_7$ of size divisible by
$
\text{LCM}\bigl( \binom{7}{2}, \binom{7}{3} \bigr)  = 105
$.
By the classification of transitive subgroups of $S_7$ \cite[p.~60]{dixon}
it follows that this Galois groups contains
$A_7$.
\end{proof}

\vsp

\subsection{Interpretation in terms of fiber products}
       \label{subsec:fiber}

\mbox{}

We continue with the notation of the previous subsection and assume properties
(i)-(iii) are satisfied.  Fix a labelling
$\lambda_1, \ldots, \lambda_n$ of the roots of $\lambda=\Lambda_1$
in $K'$, and let $\Sigma=\Sigma_1=\{\lambda_1, 
\ldots, \lambda_n\}.
$
For an integer
$
j\in [1, n-1]
$, let $\Sigma_j$ be the set of roots of $\Lambda_j$ in $K'$, and let
$\Sigma^{(j)}$ be the set of ``$j$-subsets'' of $\Sigma$ (i.e. those of
cardinality $j$).  Recall that $\Lambda_j$ splits into linear factors
over $K'$, hence $\# \Sigma_j=\# \Sigma^{(j)}=\binom{n}{j}$.
Each of these sets carries a natural action of $\gal(K'/K_0)\simeq S_n$.

\vsp

\begin{lem}
	\label{lem:obvious}
For each $j\in [1,n-1]$, there is a bijective correspondence between
$\Sigma_j$ and $\Sigma^{(j)}$ which respects the natural action of 
$\gal(K'/K_0)$ on these sets.
\end{lem}

Before proving Lemma \ref{lem:obvious}, let us state two applications of
it that we shall need.

\begin{prop}
        \label{prop:obvious}
For $\alpha \in K$, the $K$-rational roots of $\Lambda_j(x,\alpha)$ 
are in one-to-one correspondence
with the $K$-rational degree $j$ factors of $\lambda(x,\alpha)$.
\end{prop}
\begin{proof}  The $K$-rational linear factors of $\Lambda_j(x,\alpha)$
are in one-to-one correspondence with the fixed points of $G_\lambda$
in its action on $\Sigma^{(j)}$.  By the Lemma, these are in one-
to-one correspondence with the $G_\lambda$-invariant subsets of $\Sigma$
of size $j$.   The roots of a $K$-rational degree $j$ factor of 
$\lambda(x,\alpha)$ clearly form such a subset, and conversely,
a $G_\lambda$-invariant $T\in \Sigma_j$ gives the $K$-rational degree
$j$ factor $\prod_{\theta \in T} (x-\theta)$ of $\lambda$.
\end{proof} 

\begin{remark}\label{rem:prop1}
Proposition \ref{prop:obvious} lends some perspective on Proposition
\ref{prop:muller0}.  Namely, $\lambda$ has a degree $j$ factor, $1\leq
j \leq n-1$, over some intermediate field $K_0 \subseteq E \subseteq
K'$ if and only if $\Lambda_j$ has a root in $E$, i.e. if and only if
$E$ contains (a conjugate of) $K_j$.  Thus the hypothesis of Proposition
\ref{prop:muller0}, namely that $g(X_E)\geq 2$ for every $E$ over which
$\lambda$ is reducible is equivalent to the hypothesis that $g(X_j)\geq 2$
for $1\leq j \leq n-1$.  One then obtains Proposition \ref{prop:muller0}
by applying Proposition \ref{prop:obvious} in conjunction with Faltings'
Theorem.
\end{remark}

\begin{prop}
	\label{prop:an}
Suppose $1\le j\leq n-1$.  Then $\Lambda_j(x, t)$ is irreducible over
the subfield of $K'/K_0$ fixed by $A_n$.
\end{prop}

\begin{proof}[Proof of Proposition \ref{prop:an}]
Since $A_n$ is $(n-2)$-transitive, if $2\leq j \leq n-2$ 
then $A_n$, as a subgroup of
the group of permutations on the set
$
\Sigma
$,
acts transitively on the set $\Sigma^{(j)}$.  Thanks to Lemma
\ref{lem:obvious}, $A_n$, as a subgroup of
$
\gal(K'/K_0)
$,
then acts transitively on the set of roots of $\Lambda_j$ in $K'$, establishing
the Proposition for this range of $j$.

\vsp

Write $F$ for the fixed field of $K'/K_0$ by $A_n$.  If $\Lambda_1$
is reducible over $F$, then
$
\gal(K'/F)
$
is contained in
$
S_l\times S_{n-l}
$
for some $1\le l\le n-1$.  Since $F/K_0$ is quadratic, 
$
\#\gal(K'/K_0) \le 2\cdot \#S_l\cdot \#S_{n-l} < \#S_n
$, 
a contradiction.  Thus $\Lambda_1$ is irreducible over $F$.  Thanks
to Lemma \ref{lem:obvious}, that means $A_n$, as a subgroup of the group of
permutations of
$
\Sigma
$,
acts transitively on $\Sigma^{(1)}$, hence also on $\Sigma^{(n-1)}$.
Applying Lemma \ref{lem:obvious} again, 
we see that $\Lambda_{n-1}$ is irreducible
over $F$, as desired.
\end{proof}

\vsp

We now verify Lemma \ref{lem:obvious} via a fiber product
construction.  The Lemma  and the construction are
probably well-known, but we cannot locate a reference for either one so we
give the details here.  We begin with a general setup.  Recall that $K$
is a field of characteristic $0$.

\vsp

Let
$
\wp_K^n
$
denote the set of equivalence classes of non-zero, 
degree $\le n$ polynomials in
$K[x]$, where two polynomials are identified if they are $K^\times$-multiples
of each other.   We have a natural bijection between
$
\wp_K^n
$
and the set of $K$-rational points
$
\mathbf P_K^n(K)
$
of projective $n$-space, via
$$
a_0 x^n + a_1 x^{n-1} +  \cdots + a_n  \longmapsto [a_0: \cdots : a_n].
$$
In light of this, to give a polynomial
$
\lambda(x, t)\in K[x, t]
$
which is non-constant and of degree $\le n$ in $x$ is to give a non-constant
$K$-morphism
$
\Lambda: \pf{1} \rarr \mathbf P_K^n
$.
Also, for every $1\le j<n$ the multiplication
map
$
\wp_K^j\times \wp_K^{n-j}\rarr \wp_K^n
$
gives rise to a $K$-morphism
$
\phi_{n,j}: \mathbf P_K^j\times \mathbf P_K^{n-j} \rarr \pn
$,
whence a pull-back diagram
\begin{equation}
\begin{diagram}
\divide\dgARROWLENGTH by2
\node{X_j^\circ}
	\arrow[2]{e,t}{\Lambda(j)}
	\arrow{s,l}{\tilde{\phi}_{n,j}^\circ}
\node[2]{\mathbf P^j_K\times \mathbf P^{n-j}_K}
	\arrow{s,r}{\phi_{n, j}}
\\
\node{\mathbf P_K^1}
	\arrow[2]{e,t}{\Lambda}
\node[2]{\pn}
\end{diagram}
	\label{pull}.
\end{equation}
Denote by
$
\phi_n: (\mathbf P_K^1)^n \rarr \mathbf P_K^n
$
the $K$-morphism corresponding to the $n$-fold multiplication map
$
(\wp_K^1)^n \rarr \wp_K^n
$.
Then we have an analogous pull-back diagram
\begin{equation}
\begin{diagram}
\divide\dgARROWLENGTH by2
\node{\ov{X}}
	\arrow[2]{e,t}{\ov{\Lambda}}
	\arrow{s,l}{\ov{\phi}_n}
\node[2]{(\mathbf P_K^1)^n}
	\arrow{s,r}{\phi_n}
\\
\node{\mathbf P_K^1}
	\arrow[2]{e,t}{\Lambda}
\node[2]{\pn}
\end{diagram}
	\label{pkn}.
\end{equation}
Any permutation of the $n$-coordinates of the points of 
$
(\mathbf P_K^1)^n
$
is a $K$-morphism which is compatible with $\phi_n$.  Clearly
$
\deg \phi_n = n!
$,
so $\phi_n$ is a \textit{regular} branched cover with deck transformation
group $S_n$.

\vsp

Suppose $\Lambda$ corresponds to a separable, degree $n$
polynomial $\lambda$ over
$
K_0 = K(t)
$.
Then the fiber of
$
\ov{\phi}_n
$
over the generic point of
$
\mathbf P_K^1
$
consists of $n!$ pairwise distinct, \textit{ordered} $n$-tuples of the roots
of $\lambda$ over
$K_0$.  Every element of the Galois group $G_\lambda$ of $\lambda$ over $K_0$
permutes these $n$-tuples, and such a permutation gives rise to a permutation
on
$
(\mathbf P_K^1)^n
$
making the diagram (\ref{pkn}) commute.  Having fixed a labelling of the roots
of $\lambda$, we see that $G_\lambda$ is canonically identified with a
\textit{subgroup} of $S_n$.
  Since $\lambda$ is separable over $K_0$, these $n!$ $n$-tuples are
pairwise distinct, whence the scheme $\ov{X}$ is reduced.
Also,
\begin{eqnarray}
\text{$\ov{X}$ is $K$-reducible}
&
\Leftrightarrow
&
\text{the generic fiber of $\ov{\phi}_n$ is the disjoint union}
           \nonumber
\\
&&
\text{of non-trivial, $G_\lambda$-stable subsets}
           \nonumber
\\
&
\Leftrightarrow
&
\text{$G_\lambda$ does not act transitively on this fiber}
           \nonumber
\\
&
\Leftrightarrow
&
G_\lambda\subsetneq S_n.       \label{neq}
\end{eqnarray}
The multiplication map
$
(\wp_K^1)^n \rarr \wp_K^n
$
naturally factors through every
$
\wp_K^j\times \wp_K^{n-j}
$.
That means
$
\ov{\phi}_n
$
factors through
$
\phi_{n, j}
$
for every $j$; therefore the diagram  (\ref{pkn}) factors through the diagram
(\ref{pull}) for every $j$, and
$
\Phi_{n, j}
$
is also a regular branched cover with deck transformation group isomorphic to
$
S_j\times S_{n-j}
$:
$$
\xymatrix{
\ov{X} \ar[d]^{\tilde{\Phi}_{n, j}}
	\ar@/_2pc/[dd]_{\ov{\phi}_n}
	\ar[r]^{\ov{\Lambda}\mbox{\hspace{20pt}}}
&
(\mathbf P_K^1)^n
  \ar[d]^{\Phi_{n, j}}
  \ar@/^3pc/[dd]^{\phi_n}
\\
X_j^\circ
	\ar[d]^{\tilde{\phi}_{n,j}^\circ}
	\ar[r]^{\Lambda(j)\mbox{\hspace{20pt} }}
&
\mathbf P_K^j\times \mathbf P_K^{n-j}      \ar[d]^{\phi_{n,j}}
\\
\mathbf P^1_K
  \ar[r]^\Lambda
&
\mathbf P^n_K
}
$$
Finally, suppose $\lambda$ satisfies hypotheses (i)-(iii).  Then
$
G_\lambda = S_n
$,
whence the deck transformation group of
$
\ov{\phi}_n
$
is also $S_n$.  By (\ref{neq}), the scheme $\ov{X}$ is reduced and
$K$-irreducible, and so it makes sense to speak of the function field
$
K(\ov{X})
$.
Both
$
\tilde{\Phi}_{n, j}
$
and
$
\tilde{\phi}_{n, j}
$
are surjective, so
$
X_j^\circ
$
is also $K$-irreducible, and so it makes sense to speak of the function field
$
K(X_j^\circ)
$
as well, and 
$
K(\ov{X})/K(\po)
$
is an $S_n$-extension of function fields.  We have
$
\deg \tilde{\Phi}_{n, j} = j! (n-j)!
$,
and the same argument after (\ref{pkn}) shows that the deck transformation
group of
$
\tilde{\Phi}_{n, j}
$
is isomorphic to a subgroup of, and hence is exactly, $S_j\times S_{n-j}$.

\if 3\
{
We summarize the discussion above as follow.

\vsp

\begin{lem}
	\label{lem:sum}
Suppose $\lambda$ satisfies hypotheses {\rm(i)-(iii)}.  Then the extension
$
K(\ov{X})/K(\po)
$
determined by $\ov{\phi}_n$ is Galois, and for every
$
1\le j < n
$,
the subextension 
$
K(\ov{X}) / K(X_j^\circ)
$
determined by $\tilde{\Phi}_{n, j}^\circ$ is Galois with Galois group
$
S_{n, j}
$.
	\qed
\end{lem}
}
\fi

\vsp

\begin{proof}[Proof of Lemma \ref{lem:obvious}]
First, recall the notations
$
X_j, K_j
$
etc.~introduced after Lemma \ref{lem:plane2}, and the fact that $X_1$ is
given by
$
\lambda=0
$.
For any
$
t_0\in\mathbf P_{\ov{K}}^1
$,
the $\ov{K}$-rational points on the fibers of
$
\tilde{\phi}^\circ_{n, 1}
$
are in bijective correspondence with the $\ov{K}$-linear factors of
$
\lambda(x, t_0)
$,
while those on the fibers of
$
\tilde{\phi}_{n, 1}
$
are in bijective correspondence with the $\ov{K}$-rational points of the curve
$
\lambda(x, t)=0
$
with $t$-coordinates $t_0$.  These two sets are in natural bijective 
correspondence with each other; the universal property of the pullback diagram
(\ref{pull}) then implies that there is a $K$-isomorphism
$
\mu_n: X_1 \rarr X_1^\circ
$
such that
$
\tilde{\phi}_{n, 1} = \mu_n \tilde{\phi}^\circ_{n, 1}
$.
This allows us to identify the two $S_n$-extensions
$
K(\ov{X})/K_0
$
and $K'/K_0$.  The 
$
S_j\times S_{n-j}
$
subgroups in $S_n$ are pairwise conjugate, so we can identify the intermediate
subfields $K(X_j^\circ)$ with $K(X_j)$.  That means the smooth curve $X_j$ is
the canonical desingularization of $X_j^\circ$, and
$
\tilde{\phi}_{n, j}
$
is the extension of
$
\tilde{\phi}_{n, j}^\circ
$
to $X_j$, whence the
$
\gal(K'/K_0) \simeq S_n
$
action on the roots of $\Lambda_j$ over $K_0$ is the same as that on the
generic fiber of
$
\tilde{\phi}_{n, j}^\circ
$.
But the points on this generic fiber are precisely the $j$-subsets of 
$\Sigma$.
\end{proof}

\section{Generalized Laguerre Polynomials}
	\label{sec:general}

In this section we apply the machinery developed above to study
specializations of Generalized Laguerre Polynomials $ L_n^{(t)}(x) $
defined in the introduction.  In subsection \ref{subsec:geom}, we
study the singular locus of the plane curve $\mathcal{L}_n$ defined by
$L_n^{(t)}(x)=0$.  By analyzing the structure of maximal subgroups of
$S_n$, in subsection \ref{subsec:gal} we compute the genus of the
intermediate subfields of $K'/K_0$ over which $ \Lambda_j $ is
reducible.  In subsection \ref{subsec:laguerre} we combine these
ingredients to deduce Theorem \ref{thm:irr} following the strategy
outlined in sections \ref{sec:rational} and \ref{sec:spec}.

From now on, we fix $n$ and take $\lambda(x,t)=L_n^{(t)}(x)$, carrying
over all the notation ($K_0$, $K_1$, $K'$, $\oo_0$, $G_\lambda$, $B_\lambda$,
etc.) from sections \ref{sec:rational}, \ref{sec:rh}, \ref{sec:spec} to the
present setting.

\vsp
\subsection{The singular locus of $L_n^{(t)}(x)$}
    \label{subsec:geom}

\mbox{ }

Fix $n>2$.  Following Schur \cite[p.~54]{schur}, we homogenize
$
L_n^{(t)}(x)
$
by setting
\begin{eqnarray}
F_n(x, \nu, \mu)
&:=&
(-1)^n n! \mu^n L_n^{(\nu/\mu)}(x/\mu)		\nonumber
\\
&=&
x^n - \frac{k_n}{1}x^{n-1}
+
\frac{k_{n-1}k_n}{1 \cdots 2} x^{n-2}
-
\cdots
+
(-1)^n \frac{k_1 \cdots k_n}{1 \cdot 2 \cdots n},
						\label{kn}
\end{eqnarray}
where
$
k_j = j(\nu + j\mu)
$.
Let $\nl$ be the plane curve $F_n(x, \nu, \mu)=0$.  To simplify the
notation, we write
$
\partial_x F_j
$
for $\partial F_j/\partial x$.  Then we have the relations \cite[p.~54]{schur}
\begin{eqnarray}
x \partial_x F_m
&=&
m F_m + k_m F_{m-1},
	\hspace{143pt}
	(m\ge 1, F_0 := 1);	\label{deriv}
\\
F_m
&=&
(x-\nu-(2m-1)\mu) F_{m-1} - \mu k_{m-1} F_{m-2},
	\hspace{20pt}
	(m\ge 2).		\label{recur}
\end{eqnarray}
Setting $\mu=0$,  (\ref{kn}) becomes
$$
x^n - nx^{n-1}\nu + \frac{n(n-1)}{2} x^{n-2}\nu^2
-
\cdots
+
(-1)^n \nu^n
=
(x-\nu)^n.
$$
Thus $\nl $ has exactly one point along the line at infinity, namely
$
[1:1:0]
$.
Let
$
\iota_n: \nl\rarr\mathbf P_K^1
$
be the projection map defined by
$
 [x: \nu: \mu] \mapsto [\nu: \mu]
$.
%
%
%

\vsp

\begin{lem}
	\label{lem:j}
Suppose for some integer $j\in [0, n]$ and some point
$
z = [x(z): \nu(z): \mu(z)] \in\mathbf P^2_\C
$
with
$
x(z)\mu(z)\not=0
$,
we have
\begin{equation}
F_{n-j}\bigr|_z
=
\partial_x F_{n-j}\bigr|_z
= 0
\:\:
\text{ and }
\:\:
k_{n-j}\not=0.		\label{j}
\end{equation}
Then
$
F_{n-j-1}\bigr|_z = 0
$
and
$
k_{n-j-1}\not=0
$.
Moreover, if $j\le n-2$, then $\partial_x F_{n-j-1}\bigr|_z = 0$.
\end{lem}

\begin{proof}
Since $\mu(z)\not=0$, without loss of generality we can set $\mu(z)=1$.

Suppose $n\ge j+1$; then substitute into (\ref{deriv}) the first two relations
in (\ref{j}), we get
$
0 = k_{n-j}F_{n-j-1}\bigr|_z
$,
whence
\begin{equation}
F_{n-j-1}\bigr|_z=0.	\label{n1}
\end{equation}
Next, suppose $k_{n-j-1}=0$.  When we use the expansion (\ref{kn}) 
to evaluate
(\ref{n1}), we see that $x(z)=0$, a contradiction.  
Finally, suppose $n\ge j+2$.  Substituting (\ref{n1}) along with
the first relation in (\ref{j}) into (\ref{recur}), we get
$$
0 = -\mu(z) k_{n-j-1} F_{n-j-2}\bigr|_z.
$$
Substitute this and (\ref{n1}) back into (\ref{deriv}) and we get
$
x\partial_x F_{n-j-1}\bigr|_z = 0
$.
As $x(z)\not=0$, that means
$
\partial_x F_{n-j-1}\bigr|_z = 0
$.
This completes the proof of the Lemma.
\end{proof}

\vsp

\begin{lem}
            \label{lem:x1}
For $n\ge 3$ the curve $\nl$ has no finite singular point.
\end{lem}

\begin{proof}
Using the relations (\ref{deriv}) and (\ref{recur}), Schur \cite[p.~54]{schur}
showed that $F_n$, viewed as a polynomial in $x$, has discriminant
\begin{equation}
\mu^{\frac{n(n-1)}{2}} n! k_2 k_3^2 \cdots k_n^{n-1}.
	\label{disc}
\end{equation}
We are interested in the finite points on $\nl$, so for the rest of the proof
we can
set $\mu=1$.  Clearly it suffices to consider only the points on $\nl$ lying
above the branch locus of $\iota_n$.

\vsp

Suppose $z=(x_0, \nu_0)$ is a finite singular point.  By (\ref{disc}) we have
$
\nu_0\in\{-2, \ldots, -n\}
$,
and
\begin{equation}
F_n \bigr|_z
=
\partial_x F_n \bigr|_z
=
\partial_\nu F_n \bigr|_z
=
0.
	\label{sing1}
\end{equation}
We claim that $x_0\not=0$.  Suppose otherwise; set
$
\partial_\nu F_n = 0
$
and then substitute $x=0$ (recall that $\mu=1$), to get
\begin{eqnarray*}
0
=
(-1)^n \frac{\partial}{\partial \nu} \prod_{k=1}^n (\nu+k)
=
(-1)^n \sum_{m=1}^n \prod_{\substack{k=1\\k\not=m}}^n (\nu+k).
\end{eqnarray*}
Set $\nu=\nu_0$ and this becomes
$$
\prod_{\substack{k=1\\k\not=-\nu_0}}^n (\nu_0+k)=0,
$$
a contradiction.   Thus $x_0\not=0$.  
Also, if $k_n=0$, then from (\ref{kn}) we get $x_0=0$, a contradiction.
Thus $k_n\not=0$, i.e. $\nu_0 \neq -n$.
That means the hypotheses of Lemma \ref{lem:j} are satisfied for $j=0$.
Applying the lemma, we find the conditions of the lemma hold for $j=1$
as well as $\nu_0 \neq 1-n$.  Repeating this procedure, we find 
$
\nu_0 \not \in \{ -2, \cdots, -n\}
$, a contradiction.
Thus $\nl$ has no finite singular point.
\end{proof}

\begin{lem}
	\label{lem:quad}
Suppose $n\ge 2$.  Then
$
K(\sqrt{\text{disc } L_n^{(t)}(x)})
$
is a quadratic extension of $K_0$ corresponding to a smooth curve of genus
$
\bigl[ \frac{n-2}{4} \bigr]
$.
\end{lem}

\begin{proof}
Since $n\ge 2$, (\ref{disc}) says that 
$
\text{disc } L_n^{(t)}(x)
$
is a polynomial in $t$ whose square-free part has degree
$
\bigl[ \frac{n}{2} \bigr]
$,
and the Lemma follows.
\end{proof}

\if 3\
{
\vsp

\begin{lem}
  \label{lem:absirr}
$\nl$ is absolutely irreducible.
\end{lem}

\begin{proof}
Suppose otherwise.  Then for every $t_0\in\Q$ and for every point $z$ in the
fiber
$
\tilde{\phi}_{n, 1}^{-1}(t_0)
$,
the extension $\Q(z)/\Q$ must have degree $<n$.  This contradicts Schur's
result \cite{schur1} that $L_n^{(0)}(x)$ is $\Q$-irreducible.
\end{proof}
}
\fi

\vsp

Recall that the notation of section \ref{sec:rh}, such as 
$
\oo_0, \oo'
$
etc.~ now applies to the case
$
\lambda(x, t) = L_n^{(t)}(x)
$.
For
$
\nu\in \be=\{-2, \ldots, -n \}\subset\po
$,
denote by $\m_\nu$ the corresponding maximal ideal in $\oo_0$.  Denote by
$\oo_1$ the coordinate ring of an affine open set of $X_1$ containing all
places lying above every $\nu$ with respect to the projection map
$
\iota_n
$.
Then (\ref{disc}) says that the restriction of $\iota_n$ to $\oo_0$ is
unramified outside
the
         $\m_\nu$,
and
Lemma \ref{lem:x1} says that the inclusion map
$
\oo_0\subset \oo'
$
is an integral extension of Dedekind domains when localized at these
$\m_\nu$.  
From (\ref{kn}) and (\ref{disc}), we see that 
$
\oo_1
$
has exactly one ramified maximal ideal lying above $\m_\nu$:
\begin{equation}
\m_\nu \oo_1 = \n_0^{|\nu|} \n_1 \cdots \n_{s},
         \label{cons}
\end{equation}
where the $\n_i$ are pairwise distinct; in other words each branch
point $\nu$ of $L_n^{(t)}(x)$ is simple of index $|\nu|$.  
Applying Lemma \ref{lem:ram}, we
deduce the following
result.

\vsp

\begin{lem}
	\label{lem:ev}
For $\nu \in \{ -2, \ldots, -n\}$, 
let $\M_\nu\subset\oo'$ be a maximal ideal lying above $\m_\nu$.  Then
the inertia group $I(\M_\nu/\m_\nu)$ is generated by a 
cycle of length $|\nu|$.  In particular,
$$
e_\nu := e(\M_\nu/\m_\nu) = |\nu|.
$$
	\qed
\end{lem}

\vsp

\begin{prop}
       \label{prop:gx1}
Suppose $n\ge 6$.  Then the geometric genus of $\nl$ is $>1$.
\end{prop}

\begin{proof}
First, assume $n\geq 7$.
Thanks to Lemma \ref{lem:ev}, we can apply Lemma \ref{r-h}
with $V=\{ -n, 1-n, \ldots, 5-n \}$.  Since there is a unique
prime in $\oo_{1}$ above $\m_\nu$ with non-trivial ramification
index $-nu$, we have
$
c_1(\nu) = n + \nu
$,
and (\ref{cd1}) becomes
\begin{equation}
g(\nl)=g(K_1)
\ge
1 + \frac{n}{2}\Bigl( -2 + \sum_{i=0}^5 \Bigl(1-\frac{1}{d(n-i)}\Bigr) \Bigr)
-
\frac{1}{2}
\sum_{i=0}^5   \Bigl(1-\frac{1}{d(n-i)}\Bigr)\cdot i
         \label{77}
\end{equation}
For any six consecutive, positive integers, exactly two of them are prime to
$6$, another one is odd, and the remaining three are even.  Thus the first
$i$-sum in (\ref{77}) is
$$
\ge -2 + 6 - 2\times\frac{1}{5} - \frac{1}{3} - 3\times \frac{1}{2}
=
\frac{53}{30}.
$$
Thus (\ref{77}) yields
\begin{equation}
g(\nl)
\ge
1 + \frac{53n}{60} - \frac{1}{2}\Bigl(1-\frac{1}{n}\Bigr)(0+1+2+3+4+5),
         \label{777}
\end{equation}
which is $>1$ if $n>10$.  Using the more refined version (\ref{77}) we find
that in fact
$
g(\nl)>1
$
if $n\ge 7$.   Using the full Riemann-Hurwitz formula, or
the {\sc Algcurves} package in {\sc Maple}, we find that
$
g(\mathcal L_6)=4
$.
This completes the proof of the Proposition.
\end{proof}

\begin{remark}
By analyzing the singularity at infinity, one 
in fact has a nice formula $g(\mathcal L_n)=[ (n-2)^2/4]$
valid for all $n$.
\end{remark}

\vsp

\subsection{Genus of maximal subgroups}
       \label{subsec:gal}

\mbox{ }

In this subsection, we carry out the calculations which will be
necessary ingredients for the application of Theorem \ref{thm:special}
to $L_n^{(t)}(x)$ in the next section.  This involves a mixed strategy
in the following sense.  For $n\geq 10$, we show that every {\em minimal}
intermediate subfield $E$ of $K'/K_0$ has genus $>1$, thanks to Lemma
14 and Proposition \ref{prop:max} below.  It then follows from Riemann-Hurwitz
that {\em every} proper intermediate subfield has genus $>1$.
For $6\leq n \leq 9$, the
quadratic extension inside $K'/K_0$ has genus $0$ or $1$, but we have
shown in Proposition \ref{prop:an} that $\Lambda_j$ is not reducible
over this field.  It remains, then, to check for $6\leq n \leq 9$
that {\em proper} subgroups of $A_n$ give fixed fields of genus $>1$,
and this is the content of Proposition \ref{prop:alt}.  We treat $n=5$
``by hand.''


\begin{prop}
     \label{prop:max} Suppose $n\geq
6$.  If $\E$ is a maximal subgroup of $G_\lambda$ other than $A_n$, 
with corresponding fixed field $E$, then $g(X_E)>1$.
\end{prop}


\begin{prop}
   \label{prop:alt}
Suppose $6\le n\le 9$.  If $\E$ is a proper 
maximal subgroup of $A_n \subset G_\lambda$, 
with corresponding fixed field $E$, then $g(X_E)>1$.
\end{prop}

\if 3\
{
  When
we combine
Proposition \ref{prop:an} with the remark at the end of section \ref{sec:rh},
we see that it suffices to examine (\ref{cd1}) as $\E$ runs through the maximal
subgroups
of $S_n$ not containing $A_n$. 
}
\fi

\vsp

\begin{proof}[Proof of Proposition \ref{prop:max}]
Up to conjugation, the maximal subgroups of $S_n$ other than $A_n$ belong to
exactly one of the following three types
	\cite[p.~268]{dixon}:
\begin{itemize}
\item
imprimitive  subgroups: the wreath products $S_j\wr S_{n/j}$ in its
	imprimitive
	action\footnote{i.e.~the stabilizer of a partition of $n$ letters
	                 into $n/j$ disjoint subsets of equal size},
        for
	some divisor $j$ of $n$, $1<j<n$;
\item
intransitive subgroups: $S_{n, j}$ for some $1\le j<n/2$  (note that if $n$
        is even
	then $S_{n, n/2}$ is contained in $S_{n/2}\wr S_2$);
\item
a primitive subgroup of $S_n$.
\end{itemize}
For each of the three types of $\E$, we use group-theoretic properties of
$\E$ plus ramification data of
$
K'/K_0
$
to bound (\ref{cd1}) from below for large $n$, and then handle the remaining
cases individually.   Note that among any four consecutive integers $\ge 2$,
exactly one of them is prime to $6$, another one is odd, and the other two
are even.  Recall the notation $d(e)$ from Definition \ref{def1} and we see
that for $n\ge 6$,
\begin{equation}
- 2
+
\sum_{j=0}^3 \Bigl(1 - \frac{1}{d(n-j)}\Bigr) 
\ge
-2 + 4 - \frac{1}{2} - \frac{1}{2} - \frac{1}{3} - \frac{1}{5}
=
\frac{7}{15}.		\label{consec}
\end{equation}
We will also make repeated use of the
following remark.  For the rest of this section we will take
$$
V=\{3-n, 2-n, 1-n, -n\},
$$
so if $n\ge 6$ then Lemma \ref{lem:ev} implies that the inertia group of
any
$
\nu\in V
$
is generated by a single \textit{cycle}, which will allow us to use
Lemma \ref{lem:c1bound} in conjunction with (\ref{cd1}).

\vsp

\noindent
\framebox{{\bf Case: } imprimitive subgroups $S_j\wr S_{n/j}$}

First, suppose $n\ge 7$.  Since $n-3>n/2$, $S_j\wr S_{n/j}$ does not
contain any $(n-\mu)$-\textit{cycle} for $0\le \mu\le 3$.  That means
$
c_1(\nu)=0
$
for every $\nu\in V$.  Recall (\ref{consec}) and (\ref{cd1}) becomes
$$
g(X_E) \ge 1 + \frac{7}{30}[S_n: S_j\wr S_{n/j}]  > 1.
$$

\vsp

Next, suppose $n=6$.  The same reasoning as above shows that
$
c_1(\nu)=0
$
if $\nu\le 2-n$, and if $j=2$, then $c_1(3-n)=0$ as well.  So as before
$g(X_E)>1$ if
$
\E \simeq S_2\wr S_3
$.
It remains to consider the case 
$
\E \simeq S_3\wr S_2 \simeq (S_3\times S_3)\rtimes \Z/2
$.
A representative of the non-trivial coset of
$
S_3\times S_3
$
in
$
(S_3\times S_3)\rtimes \Z/2
$
(as a subgroup of $S_6$) is $(14)(25)(36)$; from this we check that elements
in this non-trivial coset all have even order.  Thus the order $3$ elements in 
$
(S_3\times S_3)\rtimes \Z/2
$
are all contained in $S_3\times S_3$.  The latter has a unique Sylow
$3$-subgroup, namely
$
\Z/3\times \Z/3
$,
so $\E$ has four distinct $\Z/3$-subgroups, whence (\ref{c11}) gives
$
c_1(-3) = \frac{8}{\#S_3\wr S_2}3\cdot 3! = 2
$.
Thus
$$
g(X_E)
\ge
1 + \frac{6!/72}{2}\frac{7}{15}
	-
	\frac{1}{2}\cdot 2 \cdot \Bigl( 1 - \frac{1}{3} \Bigr)
>1,
$$
as desired.

\vsp

\noindent
\framebox{{\bf Case: } intransitive subgroups $S_{n, j}$ with $1\le j<n/2$}

For $j>3$, $S_{n, j}$ contains no cycle of length $\ge n-3$,
so $c_1(\nu)=0$ for every $\nu\in V$.  Thus (\ref{cd1}) gives
$
g(X_E)>1
$.

\vsp

Next, suppose $j=3$, so that we can take $n\geq 7$.  
Then $c_1(\nu)=0$ for $|\nu|>n-3$, and (\ref{c12}) gives
$
c_1(3-n) < 6(n-3)
$.
Thus (\ref{cd1}) becomes
$$
g(X_E)
\ge
1 + \frac{7}{30} \frac{n!}{3!(n-3)!}
-
\frac{6n-19}{2}
\Bigl(
  1 - \frac{1}{d(n-3)}
\Bigr),
$$
which is easily 
seen to be $>1$ for $n\ge 7$ (for $n\geq 8$, use the trivial bound
$d(n-3)\leq n-3$).  

\vsp

Now, take $j=2$.  Since $n\ge 6$, the only cycles of order $n-2$ and
$n-3$ in
$
S_{n, 2} = S_2\times S_{n-2}
$
come from the cycles in $S_{n-2}$ of such order.
There are $(n-3)!$ and $(n-2)(n-4)!$ of them, respectively, so by
(\ref{c11}),
$$
c_1(n-2) = 1
\:\:
\text{ and }
\:\:
c_1(n-3) = 3,
$$
whence (\ref{cd1}) plus (\ref{consec}) gives
$$
g(X_E) \ge 1 + \frac{7}{30}\frac{n(n-1)}{2}
-\frac{1}{2}
\Bigl( 1 - \frac{1}{d(n-2)}
\Bigr)
-
\frac{3}{2} \Bigl( 1 - \frac{1}{d(n-3)}
\Bigr).
$$
This is $>1$ for $n\ge 5$, so we are done.

\vsp

Finally, consider the case $j=1$.   Then $X_E$ is simply the curve $X_1$,
which we saw right before the statement of Lemma \ref{lem:muller2} is simply
the curve $\nl$ defined by 
$
L_n^{(t)}(x)
$.
By Proposition \ref{prop:gx1}, this curve has geometric genus $>1$ if $n\ge 6$,
so we are done.

\vsp

\noindent
\framebox{{\bf Case:} primitive subgroups}

Let $\E\subset S_n$ be a primitive subgroup other than $A_n$.
By Bochert's theorem \cite[p.~79]{dixon},
$$
[S_n: \E]
\ge
\Bigl[ \frac{n+1}{2}  \Bigr]!.
$$
Using (\ref{c12}) together with the trivial estimate
$
1-\frac{1}{d(e_\nu)} \le 1 - \frac{1}{n},
$
(\ref{cd1}) becomes
\if 3\
{
Liebeck and Pyber \cite[Thm.~2]{liebeck} show that any subgroup of $S_n$
has
$
\le 2^{n-1}
$
conjugacy classes of elements.  Lemma \ref{lem:ev} implies that if $n\ge 9$ and
$
\nu, \nu'\in V = \{7-n, 6-n, \ldots, -n\}
$
are distinct then $e_\nu\not=e_{\nu'}$.  Consequently,
\begin{equation}
2^{n-1} \ge c_1(-n) + \cdots + c_1(6-n) + c_1(7-n)
         \hspace{20pt}
	 \text{if $n\ge 9$.}         \label{great}
\end{equation}
For any eight consecutive, non-zero integers, at least two of them are prime
to $6$, and at least two others are odd.  Trivially $d(e_\nu) \le n$.
Combine everything and (\ref{cd1}) becomes
}
\fi
\begin{eqnarray}
g(X_E)
&\ge&
1 + \frac{7}{30}[S_n : \E]
  -
  \frac{1}{2}\Bigl( 1 - \frac{1}{n} \Bigr)((n-1) + (n-2) + (2n-5)) + (6n-19) )
                                                         \nonumber
\\
&\ge&
1 + \frac{7}{30}\Bigl[\frac{n+1}{2}\Bigr]!
  -
  \Bigl( 1 - \frac{1}{n} \Bigr)\frac{10n-27}{2}          \label{cd21}
  \\
&\ge&
1 + \frac{7\sqrt{\pi n}}{30}\Bigl( \frac{n}{2e} \Bigr)^{n/2}
  -
  \Bigl( 1 - \frac{1}{n} \Bigr)\frac{10n-27}{2}
	\hspace{20pt}
	\text{Stirling formula \cite[p.~24]{artin}}.     \label{cd22}
\end{eqnarray}  From
(\ref{cd22}) we get that $g(X_E)>1$ if $n\ge 15$.  Using the sharper
form (\ref{cd21}), we see that in fact $g(X_E)>1$ if $n\ge 11$.  For
$n=9,10$, if we use the original inequality (\ref{cd1}), we also
obtain $g(X_E)>1$.
To handle the remaining values of $n$, i.e. $6,7,8$, 
we make use of classification of
primitive groups of small degree \cite{butler}.

\vsp

\if 3\
{
\framebox{$n=10$}
\hspace{5pt}  
The maximal primitive subgroups of $S_{10}$ other than
$A_{10}$ are
$
P\Gamma L_2(\F_9)
$,
$
S_6
$,
$
M_{10}
$,
and
$
PGL_2(\F_9)
$.
The first one has order $1440$, and the other three, $720$, hence
$
[S_{10}: \E] \ge 10! / 1440
$.
Substitute this into (\ref{cd20}) and we find that
$
g(X_E) > 1
$.
}
\fi

\vsp

\framebox{$n=8$}  
\hspace{5pt}  
$S_8$ has two maximal primitive subgroups other than
$A_7$, namely
$
PGL(2, \F_7)
$
and
$
2^3\cdot PSL_2(\F_7)
$
(a group with normal subgroup $(\Z/2)^3$ and with quotient $PSL_2(\F_7)$).
In particular, both groups contain no element of order $5$, so the $c_1$-term
in (\ref{cd1}) corresponding to the branched point $\nu=-5$ is zero.  For the
group $PGL_2(\F_7)$, (\ref{cd1}) then becomes
$$
1 + \frac{7}{30}\frac{8!}{336}
-
\frac{1}{2}
\Bigl(
  7\Bigl(1-\frac{1}{2}\Bigr)
  +
  6\Bigl(1-\frac{1}{7}\Bigr) + 11\Bigl(1-\frac{1}{2}\Bigr)
\Bigr) > 1.
$$
To handle the group $2^3\cdot PSL_2(\F_7)$ we need to refine our estimate for
the $c_1(-7)$-term.  Sylow theory dictates that
$
2^3 \cdot PSL_2(\F_7)
$
has at most $64$ Sylow $7$-subgroups, all of order $7$, so 
$
2^3 \cdot PSL_2(\F_7)
$
has at most $64\times 6 = 384$ elements of order $7$.  Substitute this into
(\ref{c11}) and we find that
$
c_1(-7) \le 2
$,
whence (\ref{cd1}) becomes
$$
1 + \frac{7}{30}\frac{8!}{8\times 168}
-
\frac{1}{2}
\Bigl(
  7\Bigl(1-\frac{1}{2}\Bigr)
  +
  2\Bigl(1-\frac{1}{7}\Bigr) + 11\Bigl(1-\frac{1}{2}\Bigr)
\Bigr) > 1.
$$

\if 3\
{
The first one has order
$168$ and $6$ conjugacy classes, and the second one, order $1344$ and $11$
conjugacy classes.
Set 
$
V=\{ 5-n, 4-n, \ldots, -n\}
$;
then the analogue of (\ref{great}) becomes
$
11 \ge c_1(5-n) + c_1(4-n) + \ldots + c_1(-n)
$,
hence (\ref{cd1}) gives
$$
g(X_E)
\ge
1
+
\frac{1}{2} \Bigl( - 2 + \sum_{j=0}^5 \Bigl( 1 - \frac{1}{d(8-j)}\Bigr)  \Bigr)
\frac{8!}{1344}
-
11 > 1.
$$
}
\fi

\vsp

\framebox{$n=7$}
\hspace{5pt}  
  $S_7$ has a unique maximal primitive subgroup other than
$A_7$, namely $PSL_2(\F_7)$.  It has $42$ elements of order $4$, no element of
order $5$, and $48$ elements of order $7$, so
$
c_1(-4) = \frac{42}{168} 4\cdot 6 = 6
$,
$c_1(-5)=0$, and
$
c_1(-7) = \frac{48}{168} \cdot 7 = 2
$,
whence (\ref{cd1}) becomes
$$
g(X_E)
\ge
1 + \frac{7}{30}\frac{7!}{168}
-
\frac{1}{2}
\Bigl(
  2(1-\frac{1}{7})
  +
  5\Bigl(1-\frac{1}{2}\Bigr)
  +
  6\Bigl(1-\frac{1}{2}\Bigr)
\Bigr)
> 1.
$$

\vsp

\framebox{$n=6$}
\hspace{5pt}  
  $S_6$ has a unique maximal primitive subgroup other than
$A_6$, namely
$
PGL_2(\F_5) \simeq S_5 \simeq S_{6, 1}
$.
For such intransitive groups we already saw that $g(X_E)>1$, so we are done
for $n=6$.
This completes the proof of Proposition \ref{prop:max}.
\end{proof}

\if 3\
{
?????????

It has $24$ elements of order $5$, 

how many elements of order $6, 5, 4, 3$ are there in $PGL_2(\F_5)$??

$$
g(X_E)
\ge
1 + \frac{7}{30}[S_6: PGL_2(\F_5)]
-
\frac{1}{2}
    \Bigl(
       ???
    \Bigr)
\ge
1 + \frac{7}{5} - ??? 
$$
as desired.  \framebox{\tt  not good.  might have to actually count
double cosets!!!!}
}
\fi

\vsp

\begin{proof}[Proof of Proposition \ref{prop:alt}]
We will make extensive use of the \textit{Atlas} \cite{atlas}
to determine the maximal subgroups of these $A_n$, and for the number of
conjugacy classes of elements $A_n$ and $PSL_2(\F_q)$.  For the rest of the
proof we take
$
V = \{-n, 1-n, 2-n, 3-n \}
$.

\vsp

\framebox{$n=9$}
\hspace{5pt}  
According to the \textit{Atlas}, the
  maximal
  subgroups\footnote{in what follows we will consider the \textit{isomorphism}
     classes, and not \textit{conjugacy classes}, of maximal subgroups of
     these $A_n$.  For the purpose of computing $g(X_E)$ this is sufficient.}
of
$A_9$ are $A_8, S_7$,
plus others of indices $\ge 84$ in $A_9$.  First, consider those $\E$ of
index
$\ge 84$ in $A_9$.  Then
$
[S_9: \E] \ge 168
$,
and (\ref{cd1}) becomes
\begin{eqnarray*}
g(X_E)
&\ge&
1 + \frac{7}{30}168
-
\frac{1}{2}
\Bigl(
  \Bigl(1-\frac{1}{3}\Bigr)c_1(-9)
  +
  \Bigl(1-\frac{1}{2}\Bigr)c_1(-8)
  +
  \Bigl(1-\frac{1}{7}\Bigr)c_1(-7)
  +
  \Bigl(1-\frac{1}{2}\Bigr)c_1(-6)
\Bigr)
\\
&\ge&
1 + \frac{196}{5}
-
\frac{1}{2}
\Bigl(
  \frac{2}{3}8 + \frac{1}{2}7 + \frac{6}{7}13 + \frac{1}{2}35
\Bigr)
> 1,
\end{eqnarray*}
which is satisfactory.  Next, take $\E=A_8$.  Then
$
[S_9: \E] = 18
$,
and $A_8$ has no \textit{cycles} of order $9, 8$ or $6$, so
$
c_1(-9) = c_1(-8) = c_1(-6) = 0
$.
There are $8!/7$ elements of order $7$ in $A_8$,
so
$
c_1(-7) = \frac{8!/7}{8!/2}7\cdot 2 = 4
$.
Thus
$$
g(X_E)
\ge
1 + \frac{7\cdot 18}{30}
-
\frac{4}{2}\Bigl(1-\frac{1}{7}\Bigr)
> 1.
$$
Finally, take $\E=S_7$.  Then $[S_9: \E] = 72$ and $S_7$ has no element of
order $9$ or $8$, so
$$
g(X_E)
\ge
1 + \frac{7\cdot 72}{30}
-
\frac{1}{2}
\Bigl(
  \Bigl(1 - \frac{1}{7}\Bigr)13
  +
  \Bigl(1 - \frac{1}{2}\Bigr)35
\Bigr)
> 1.
$$
This completes the case $n=9$.

\vsp

\framebox{$n=8$}
\hspace{5pt}  
The maximal subgroups of $A_8$, along with their indices in $A_8$, are
$$
(A_7, 8);
\:\:
( (2^3:PSL_2(\F_7)), 15);
\:\:
(S_6, 28);
\:\:
( 2^4: (S_3\times S_3), 35);
\:\:
( (A_5\times 3): 2, 56).
$$  From
(\ref{c12}) we get the standard estimates
\begin{equation}
c_1(-8)< 8,  \:\:  c_1(-7)< 7,  \:\: c_1(-6)< 6\cdot 2!.
     \label{88}
\end{equation}

\vsp

The case
$
\E =  2^3:PSL_2(\F_7)
$
has already been dealt with in the course of proving Prop.~\ref{prop:max}.
For
$
\E = 2^4: (S_3\times S_3)
$,
it has no element of order $5$ or $7$, whence $c_1(-5)=c_1(-7)=0$.  We have
$
[S_n: \E] = 70
$,
so (\ref{cd1}) becomes
$$
g(X_E)
\ge 
1 + \frac{7}{30}70
-
\frac{1}{2}
\Bigl(
  \Bigl(1-\frac{1}{2}\Bigr) 8 + \Bigl(1-\frac{1}{2}\Bigr) 12
\Bigr)
> 1.
$$
Next, take $\E=(A_5\times 3):2$, i.e.~a split extension with kernel 
$
A_5\times\Z/3
$
and quotient $\Z/2$.  The order $5$ elements in $\E$ are all in
$
A_5\times\Z/3
$,
and hence there are $4!$ of them.  Thus (\ref{c11}) gives
$
c_1(-5) = \frac{4!}{360}5\cdot 3! = 2
$.
Also, $\E$ has no element of order $7$, so $c_1(-7)=0$.  Thus (\ref{c12})
becomes
$$
g(X_E)
\ge
1 + \frac{7}{30}112
-
\frac{1}{2}
\Bigl(
  \Bigl(1 - \frac{1}{2}\Bigr)8 + \Bigl(1 - \frac{1}{2}\Bigr)6
  +
  \Bigl(1 - \frac{1}{5}\Bigr)2
\Bigr) > 1.
$$
For $\E=S_6$, again it has no order $7$ elements so $c_1(-7)=0$.  It has
$6!/5$ order $5$ elements, so
$
c_1(-5) = \frac{6!/5}{6!}5\cdot 3! = 6
$.
Thus (\ref{cd1}) becomes
$$
1 + \frac{7}{30}56
-
\frac{1}{2}
\Bigl(
  \Bigl(1 - \frac{1}{2}\Bigr)8 + \Bigl(1 - \frac{1}{2}\Bigr)6
  +
  \Bigl(1 - \frac{1}{5}\Bigr)6
\Bigr) > 1.
$$
Now take $\E=A_7$.  There are no cycles of length $6$ or $8$ in $A_7$, so
$
c_1(-8)=c_1(-6)=0
$.
There are $6!$ order $7$ elements and $7!/(5\cdot 2!)$ order $5$ elements in
$A_7$, so
$
c_1(-7) = c_1(-5) = 1
$.
Thus
$$
g(X_E)
\ge
1 + \frac{7}{30}8 
-
\frac{1}{2}
\Bigl(
  \Bigl(1-\frac{1}{7}\Bigr) + \Bigl(1-\frac{1}{5}\Bigr)
\Bigr) > 1.
$$

\vsp

\framebox{$n=7$}
\hspace{5pt}  
The maximal subgroups of $A_7$, along with their indices in $A_7$, are
$$
(A_6, 7);
\:\:
(PSL_2(\F_7), 15);
\:\:
(S_5, 21);
\:\:
((A_4\times 3): 2, 35).
$$
Note that (\ref{c12}) gives the following estimates
$$
c_1(-7)<7, \:\: c_1(-6) < 6, \:\: c_1(-5) < 5\cdot 2, \:\: c_1(-4) < 4\cdot 6.
$$
First, take
$
\E = (A_4\times 3): 2
$.
Then $\E$ has no element of order $7$ or $5$, so
$
c_1(-7) = c_1(-5) = 0
$. 
Thus (\ref{cd1}) becomes
$$
g(X_E)
\ge
1 + \frac{7}{30}70
-
\frac{1}{2}
\Bigl(
  \Bigl(1-\frac{1}{2}\Bigr)5 +  \Bigl(1-\frac{1}{2}\Bigr)23
\Bigr) > 1.
$$
Next, take $\E=S_5\subset A_7$.  Then it has no \textit{cycles} of order
$7$ or $6$, so
$
c_1(-7)=c_1(-6)=0
$.
It has $4!$ elements of order $5$, and $5!/4$ elements of order $4$.  Thus
$
c_1(-5) = 2
$
and
$
c_1(-4) = 6
$.
Thus
$$
g(X_E) \ge 1 + \frac{7}{30}42 - \frac{1}{2}
\Bigl(
  \Bigl(1-\frac{1}{5}\Bigr)2 +  \Bigl(1-\frac{1}{2}\Bigr)6
\Bigr) > 1.
$$
Now, take $\E=A_6\subset A_7$.  It has no order $7$ elements and no
\textit{cycles} of order $6$ or $4$.  It has $6!/5$ order $5$ elements, so
$
c_1(-5) = 2
$.
Thus
$$
g(X_E) \ge 1 + \frac{7}{30}14 - \frac{1}{2}\Bigl(1-\frac{1}{5}\Bigr) 2 > 1.
$$
Finally, take
$
\E = PSL_2(\F_7)
$.
It has $42$ elements of order $4$, none of order $5$ or $6$, and $48$ elements
of order $7$.  Thus
$
c_1(-4) = \frac{42}{168}4\cdot 3! = 6,
c_1(-5) = c_1(-6)=0,
c_1(-7) = \frac{48}{168}7 = 2
$.
Then
$$
g(X_E) \ge 1 + \frac{7}{30}30 - \frac{1}{2}
\Bigl(
  \Bigl(1-\frac{1}{7}\Bigr)2 + \Bigl(1-\frac{1}{2}\Bigr)6
\Bigr) > 1.
$$

\vsp

\framebox{$n=6$}
\hspace{5pt}  
The maximal subgroups of $A_6$, along with their indices in $A_6$, are
$$
(A_5, 6);
\:\:
( (\Z/3\times\Z/3)\rtimes\Z/4, 10);
\:\:
(S_4, 15).
$$
First, take $\E=S_4$.  It has six elements of order $4$, eight of order $3$,
and none of order $5$ or $6$.  Thus $c_1(-4)=2, c_1(-3)=6, c_1(-6)=c_1(-5)=0$,
whence $g(X_E)>1$.

\vsp

Next, take $\E=A_5$.  It has twenty-four elements of order $5$, twenty elements
of order $3$, and none of order $6$ or $4$.  Thus
$
c_1(-5)=1, c_1(-3)=3, c_1(-6)=c_1(-4)=0
$,
whence
$
g(X_E)>1
$.

\vsp

Finally, take 
$
\E=(\Z/3\times\Z/3)\rtimes\Z/4
$.
Then $c_1(-5)=0$.  There are $8$ elements of order $3$, and hence
$
\le 27
$
elements of order $4$.  Thus
$
c_1(-3) = 4
$
and
$
c_1(-4) \le 6
$.
It follows that $g(X_E)>1$.  This completes the proof of Proposition 
\ref{prop:alt}.
\end{proof}

\vsp

\subsection{Proof of Theorem \ref{thm:irr}}
     \label{subsec:laguerre}

\mbox{ }

\framebox{\sc Step I.}
First, we treat the case $n=5$ using an argument specific to
quintics. A separable quintic over $K$ (not necessarily irreducible)
has a solvable Galois group if and only if its resolvent sextic has a
root in $K$ \cite{dummit}.  Compute the resolvent sextic of
$
L_5^{(t)}(x)
$
using the formula in \cite{dummit} and set it to
$
(x-10A)(x^5 + c_1x^4 + \cdots + c_5)
$,
obtaining six equations in $t, A, c_1, \ldots, c_5$.
Eliminate 
$
c_1, \ldots, c_5
$ from
the six equations using {\sc Maple} and we arrive at a
single equation in
$t$ and
$A$:
\begin{eqnarray*}
&&
A^6+(-12t^2-24t)A^5+(120t^2+60t^3)A^4
\\
&&
\hspace{16pt}
+(720t^3+2120t^4+1600t^5+360t^6)A^3
\\
&&
\hspace{16pt}
+(-5040t^4-11580t^6-4200t^7-540t^8-13200t^5)A^2
\\
&&
\hspace{16pt}
+(10368t^4+39744t^5+48864t^6+14448t^7-12480t^8-9360t^9-1728t^{10})A
\\
&&
\hspace{16pt}
-3(5832t^5+26892t^6+50814t^7+50645t^8+28406t^9+8735t^{10}+1278t^{11}+54t^{12}).
\end{eqnarray*}
Using the {\sc Algcurves} package in {\sc Maple}, we find that this equation
is absolutely irreducible and defines a plane curve with geometric genus $3$.
Thanks to Faltings, that means 
$
L_5^{(\alpha)}(x)
$
is $K$-irreducible and is not solvable for $\af$.
This completes the proof for the case $n=5$. From now on, assume that $n\ge 6$.

\vsp

\framebox{\sc Step II.}
Given a number field $K$, we claim that if there exists \textit{one}
$
\beta\in K
$
for which
$
L_n^{(\beta)}(x)
$
has $S_n$-Galois group over $K_0$, then Theorem \ref{thm:irr} holds for this
$K$.

\vsp

By (\ref{disc}), the discriminant of $L_n^{(t)}(x)$ is not constant.
Since $n\ge 5$, Lemma \ref{lem:alg} applies so that the existence of
this one $\beta$ yields the necessary hypotheses on $K'/K_0$.  For
$n\geq 10$, the genus of the fixed field of every proper maximal
subgroup of $G_\lambda$ is greater than $1$ (Proposition
\ref{prop:max} and Lemma \ref{lem:quad}).  By Riemann-Hurwitz, since
$K$ has characteristic $0$, $g(X_E)\leq g(X_{E'})$ whenever $E \subset
E'$.  Thus, for $n\geq 10$, {\em every} non-trivial intermediate
subfield of $K'/K_0$ has genus greater than $1$.  For degrees
$n=6,7,8,9$, we have shown, (a) that {\em proper} maximal subgroups of
$A_n$ and $S_n$ have genus greater than one (Propositions \ref{prop:max} and
\ref{prop:alt}), and (b) over the quadratic subfield of $K_0$ in $K'$, 
the polynomials $\Lambda_j$ are all irreducible (Proposition \ref{prop:an}).
Thus, for all $n\geq 6$, the hypotheses of Theorem 3 are satisfied.

We therefore obtain the first part of Theorem
\ref{thm:irr}(a) for $n\ge 7$.  By Lemma \ref{lem:quad}, if $n\ge 10$
(resp.~$n\ge 6$) then the set of $t\in K$ corresponding to even Galois groups
are parameterized by a curve of geometric genus $\ge 2$ (resp.~$\ge 1$).
The rest of Theorem \ref{thm:irr} for $n\ge 7$ now follows. 

\vsp

For $n=6$, the argument for Theorem \ref{thm:special} only shows that the
degree of the splitting field of all but finitely many
$
L_n^{(\alpha)}(x)
$
over $K$ is divisible by 
$
\text{LCM}\bigl( \binom{6}{2}, \binom{6}{3} \bigr) = 60
$.
To improve this we use a different test function.  By Lemma
\ref{lem:alg}(a), the fixed field of $K'/K_0$ by
$
S_3\times \{1\}\subset S_{6, 3}
$
corresponds to a smooth projective curve $X_{3,0}$ plus a $K$-morphism
$
\xi_{3,0}: X_{3,0}\rarr\po
$.
Write
$
\Lambda_{3, 0}(x, t) = 0
$
for the corresponding birational plane curve.  The same argument as in Lemma
\ref{lem:obvious} shows that the roots of $\Lambda_{3, 0}(t)$ over $K_0$ are
in bijective correspondence with triples of roots of
$
L_6^{(t)}(x)
$
over $K_0$.  Argue as in Proposition \ref{prop:an} and we see that
$
\Lambda_{3, 0}(x, t)
$
is irreducible over the fixed field of $K'/K_0$ by $A_6$.  The discussion in
subsection \ref{subsec:gal} is now applicable, and we see that for $\af$, the
degree of the splitting field of 
$
L_n^{(\alpha)}(x)
$
over $K$ is divisible by
$
\deg \xi_{3,0} = [S_6: S_3\times \{1\}] = 120
$.
By the classification of transitive subgroups of $S_6$ \cite[p.~60]{dixon},
we are done.

\vsp

\framebox{\sc Step III.}
Schur \cite{schur1} showed that
$
L_n^{(0)}(x)
$
is $\Q$-irreducible and has $S_n$ Galois group.  That means
$
L_n^{(t)}(x)=0
$
has $S_n$ Galois group over $\Q(t)$.  Apply Step II and we get Theorem
\ref{thm:irr} for $K=\Q$.  In particular, 
$
\lambda(x, \alpha)
$
has $S_n$ Galois group over $\Q$  for all but finitely many
$
\alpha\in\Z
$.  From
(\ref{disc}) we see that, for any finite set of primes $\Sigma$,
infinitely many of these $S_n$-extensions of $\Q$ must be ramified outside 
$
\Sigma
$.
There are only finitely many number fields of bounded degree
which are unramified outside $\Sigma$, so for any fixed number
field $K$, there exist infinitely many $\alpha'\in \Q$ so that any root of
$
L_n^{(\alpha')}(x)
$
defines a degree $n$ extension
$
L_{\alpha'}/\Q
$
with $S_n$-Galois closure and is ramified at a prime which is unramified in
$K/\Q$.  Since $S_n$ has no subgroup of index $< n$, that means
$
L_{\alpha'}\cap K= \Q
$,
whence
$
L_n^{(\alpha')}(x)
$
also has $S_n$ Galois group over $K$.  Apply Step II with $\beta=\alpha'$ and
we are done.

\vsp

\section{Simple covers}
       \label{sec:simple}

Let $Y$ be a smooth projective curve defined over a number field $K$, and let
$
\pi: Y\rarr\po
$
be a $K$-morphism of degree $n$.  We say that 
$\pi$ is a \textit{simple cover}
if the fiber above every point in
$
\mathbf P_{\ov{K}}^1
$
contains at least $n-1$ distinct points.  In other words, every branch 
point of $\pi$ is simple of index $2$.   By \cite[top of p.~549]{fulton},
the (geometric) Galois group of a simple $n$-cover is precisely $S_n$.
 Say $Y$ has genus $g$; then the
Riemann-Hurwitz formula implies that the number of branch points of
$\pi$ is exactly
\begin{equation}
\# B_\pi = 2g+2n-2.		\label{exact}
\end{equation}
Over an algebraically closed field, if $n\ge g+1$ then every smooth
projective curve of genus $g$ admits a simple cover of degree $n$
\cite[Prop.~8.1]{fulton}.

\vsp

Suppose $\lambda(x, t)\in K[x, t]$ is irreducible over $K_0=K(t)$ of
degree $n$ and defines a simple cover $K_1/K_0$ (in the notation of
section \ref{sec:rational}).  To simplify the exposition, suppose $K$
is algebraically closed in the splitting field $K'$ of $\lambda$ over
$K_0$.
The following example of M\"uller shows that we cannot expect
all but finitely $K$-specializations of $\lambda$ to be
$K$-irreducible, let alone having the same Galois group as $\lambda$.
Consider the transpositions
$
g_1 = (1, 2), g_2 = (2, 3), \ldots, g_{n-2} = (n-2, n-1), g_{n-1} = (n-1, n),
g_n = (n-1, n), g_{n+1} = (n-2, n-1), \ldots, g_{2n-3} = (2, 3), g_{2n-2} =
(1, 2)
$.
Note that the product of these $g_i$ is $1$, and that they generate $S_n$.
So by the Riemann existence theorem \cite[Cor.~7.3]{volklein}, there exists a
degree $n$ branched cover
$
X_n\rarr \mathbf P^1_{\ov{K}}
$
with exactly $2(n-1)$ branched points over $\ov{K}$, such that the inertia
group of the $i$-th branch point is generated by $g_i$.
By Riemann-Hurwitz, the cover with this description has geometric genus zero
and is a simple cover.  So taking a finite extension $L/K$ if necessary,
there are infinitely many $L$-rational specializations of this cover with an
$L$-linear factor.

\vsp

This example shows that there is not an analogue of Theorem
\ref{thm:irr} which holds for {\em all} simple covers of sufficiently
large degree.  But if we start with a simple cover of genus at least
$2$, then we can reach a similar conclusion, as in Part (b) of the
following theorem.  Even if we start with a rational or elliptic
simple cover, however, Part (a) of the theorem says that all but
finitely many specializations are either irreducible or factor as a
linear times a degree $n-1$ irreducible factor.  We give two proofs of
Theorem \ref{thm:simple}.  The first is due to M\"uller, and uses a
classification theorem of Liebeck and Saxl; we thank M\"uller for
suggesting that we include it here, as well as for catching an error
in an earlier version of the theorem.  The second proof illustrates
the usefulness of the interpretation of the curve $X_j$ introduced in
section \ref{sec:spec} as the variety whose $K$-rational points
parametrize the $K$-rational degree $j$ factors of $\lambda$.

\vsp

\begin{thm} 
          \label{thm:simple}
Let $\lambda(x, t)$ be an irreducible polynomial over
$K(t)$ defining a simple cover $\pi: Y\rarr\po$ of degree $n\geq 5$
and geometric genus $g=g_Y\geq 0$.  If $g=0$, assume $n\geq 6$.  Then,

{\rm (a)}  
For $\af$, the specialization $\lambda(x,\alpha)$ has a $K$-irreducible 
factor of degree $\geq n-1$.

{\rm (b)}
If  $g_Y\geq 2$, then for $\af$, the specialization $
\lambda(x, \alpha) $ is $K$-irreducible.   

{\rm (c)}  If $g_Y\geq 2$ and $n\geq 7$, for $\af$,
the Galois group of $\lambda(x,\alpha)$ over $K$ is $S_n$.  
\end{thm}

\begin{proof}[First Proof]
Suppose $\E$ is a subgroup of $G_\lambda$ with fixed field $E$.
The key step is the following claim.
\vsp

{\bf Claim.} If $\E$ is a maximal subgroup of $G_\lambda$ not
conjugate to $S_{n,1}$, then $g_E\geq 2$.
\vsp

We now give a proof, communicated to us by M\"uller, of this claim.
Suppose $\E$ is a maximal subgroup of $S_n$ which is 
not conjugate to $S_{n,1}$.
Recall that $K'/K_0$ is the Galois closure of the function field extension
$K_1/K_0$ defined by the simple cover $\pi$.  This yields an action of
$
\gal(K'/K_0)\simeq S_n
$
on the generic fiber of $\pi_E$.  By Galois theory, this action, call it
$\rho_E$, is simply the left-action of $S_n$ on the left cosets
of $\E$ in $S_n$.  Since
$
\E\not\simeq S_{n, 1}
$,
this action is not the natural degree $n$ action of $S_n$.  Let 
$
\mu(\E)
$
be the largest integer $m$ such that every
transposition of $S_n$ moves at least $m$ points in 
the $\rho_E$-action.  Since 
$
\pi_E
$
is a quotient of the Galois closure of the \textit{simple} cover $\pi$, the
ramification index of
$\pi_E$
at any maximal ideal $\n$ of an affine coordinate ring of $X_E$ divides $2$
(Lemma \ref{lem:ram}).  By definition of $\mu(\E)$, there are $\mu(\E)/2$
$\oo_E$-primes $\n$ above $\m$ with $e(\n/\m)=2$, 
thus for any
$
\m\in B_\pi
$,
as $\n$ runs through all maximal ideals of $\oo_E$ lying above $\m$, we have
$$
\sum_{\n/\m}  \bigl( e(\n/\m) - 1 \bigr) f(\n/\m)
\ge
\mu(\E)/2.
$$
By Lemma \ref{lem:alg}(b), the branch locus of $\pi_E$ is exactly $B_\pi$.
Thus Riemann-Hurwitz gives 
\begin{equation}
2(N_E-1+g_E) \ge \#B_\pi\times \mu(\E)/2.            \label{min}
\end{equation}

\vsp

Suppose $g_E\le 1$.  Then (\ref{min}) and (\ref{exact}) together
give 
\begin{equation}
\mu(\E) \le 2N_E/(n+g_Y-1)\leq 2N_E/(n-1)               \label{bad}
\end{equation}
Recall that $\rho_E$ is transitive, and since $\E$ is maximal,
\cite[Cor.~1.5A]{dixon} implies that $\rho_E$ is primitive as well.  
By \cite[Thm.~6.1]{liebeck-saxl}, either
\begin{itemize}
\item
$
\mu(\E)\ge N_E/2
$,
or
\item
$S_n$ contains a normal subgroup isomorphic to $H^r$, where $H$ is
isomorphic to an alternating group $A_m$ for some $m$, or to a simple group of
Lie type over $\Z/2$.
\end{itemize}
The first option plus (\ref{bad}) implies that $n\le 4$ for $g_Y\geq 1$
and $n\leq 5$ for $g_Y=0$, and we are done.
Since $n\ge 5$, for
the second option we must have
$
m=n, r=1
$
and
$
H\simeq A_n
$.
Furthermore, \cite[Thm.~6.1]{liebeck-saxl} says that, in this case, 
$\rho_E$ is in fact the
action of $S_n$ on the set of $j$-subsets of 
$
\{1, \ldots, n\}
$ for some $j\in [1,n/2]$,
whence
$
N_E=\binom{n}{j}
$,
and
$
\mu(\E)=2\binom{n-2}{j-1}
$.
Recall (\ref{bad}) and we get
\begin{equation}\label{last}
2\binom{n-2}{j-1} \le  \frac{2}{n+g_Y-1}\binom{n}{j}.
\end{equation}
This inequality simplifies to $j(n-j) \leq n(n-1)/(n+g-1)$.  
Since $n\geq 5$, this
is only possible if $g_Y\leq 1$ and $j=1$.  Thus, $g_E \geq 2$
for all maximal subgroups $\E$ of $G_\lambda$ not conjugate
to $S_{n,1}$.

\vsp

Parts (b) and (c) follow immediately from the claim.  Indeed, the hypothesis
there, namely $g_{K_1}>1$, ensures that the genus of {\em every} minimal
subfield of $K'/K_0$ is at least $1$.  Now we can apply Proposition
\ref{prop:muller0} and Theorem \ref{thm:special} to complete the proof.
For part (a), it remains only to combine the claim with 
Proposition \ref{prop:obvious}.
\end{proof}

\vsp

\begin{proof}[Second Proof]
Now we give a slightly different approach which is independent of
Liebeck-Saxl.
Fix $j\in [1,n/2]$ and suppose $\lambda(x,\alpha)$ has a $K$-rational
degree $j$ 
factor for infinitely many $\alpha \in K$. Then Proposition \ref{prop:obvious}
implies that $g(X_j)\leq 1$.  But
the function field of $X_j$ is the fixed field $E=K_j$ of
$\E=S_{n,j}$, so by Lemma \ref{lem:obvious} (recall the notation
introduced at the beginning of subsection
\ref{subsec:fiber}), we have reduced again
to the case where $\rho_E$ is the action of $G_\lambda$ on the set of
$j$-subsets of
$
\Sigma=\{\lambda_1,\ldots, \lambda_n\}
$.
Repeat the argument arising from (\ref{last}) and we get $g_Y\leq 1$ and $j=1$,
from which Parts (a) and (b) of the Theorem follow.
To prove (c), assuming now that $g_Y\geq 2$, and $n\geq 7$, 
we have already seen that the fixed field
of the {\em intransitive} 
maximal subgroups $S_{n,j}$ have genus at least 2.
Now we consider a {\em transitive} maximal subgroup $\E$ of $G_\lambda$.
By
\cite[Lem.~4.4.4]{serre:gal}, 
the only transitive subgroup of $S_n$ that contains
a $2$-cycle is $S_n$,
so we may assume $\E$ has no transpositions.  
Recalling that
$\# B_\pi = 2g+2n-2$, (\ref{cd1}) and (\ref{c11}) combine to give $g(X_E)>1$
in this case as well.  Now we can apply Theorem \ref{thm:special} to conclude
the proof.
\end{proof}

\vsp

\begin{remark}
The argument above plus Theorem \ref{thm:special} shows that 
if $g_Y\ge 2$ then the Galois group of 
$
\lambda(x, \alpha)
$
has order divisible by $60$ (if $n=6$) and by $20$ (if $n=5$) 
for $\af$.  We do not know
if the Galois groups are in fact $S_6$ and $S_5$, respectively, for $\af$.
\end{remark}

\vsp

\begin{ack}
We are grateful to Professor M\"uller for pointing out an error in the
final section of an earlier draft, and for showing us the first proof of
Theorem \ref{thm:simple}.  We would also like to thank Professors
Boston, Cox, Gunnells and Markman for useful discussions.
\end{ack}

\bibliographystyle{amsalpha}

\vfill

\end{document}